\documentclass[11pt]{article}

\usepackage{amsmath,amsthm}
\usepackage{amssymb}
\usepackage{amsfonts}
\usepackage{amscd}
\usepackage{hyperref}

\newtheorem{thm}{Theorem}[section]
\newtheorem{theorem}[thm]{Theorem}

\newtheorem{proposition}[thm]{Proposition}

\theoremstyle{definition}
\newtheorem{definition}[thm]{Definition}
\newtheorem{example}[thm]{Example}

\newtheorem{remark}[thm]{Remark}
\newtheorem{observation}[thm]{Observation}
\newtheorem{free text}[thm]{}

\newcommand{\N} {\mathbf{N}}
\newcommand{\R} {\mathbf{R}}
\newcommand{\Z} {\mathbf{Z}}
\newcommand{\C} {\mathbf{C}}
\newcommand{\T}{\mathbf{T}}
\newcommand{\cU} {\mathcal{U}}

\newcommand{\id}{\mathrm{id}}
\newcommand{\Ad}{\mathrm{Ad}}
\newcommand{\Hom}{\mathrm{Hom}}

\newcommand{\lra} {\longrightarrow}

\newcommand{\End}{\mathrm{End}}

\newcommand{\tr}{\mathrm{tr}}

\newcommand{\fsp}{\mathfrak{sp}}
\newcommand{\fso}{\mathfrak{so}}
\newcommand{\fsu}{\mathfrak{su}}

\newcommand{\fsl}{\mathfrak{sl}}

\newcommand{\fg}{\mathfrak{g}}
\newcommand{\fh}{\mathfrak{h}}

\newcommand{\ad}{\mathrm{ad}}

\newcommand{\rmt}{\mathrm{t}}

\newcommand{\al}{\alpha}

\newcommand{\ep}{\epsilon}
\newcommand{\de}{\delta}
\newcommand{\la}{\lambda}

\newcommand{\cO}{\mathcal{O}}
\newcommand{\cA}{\mathcal{A}}
\newcommand{\cB}{\mathcal{B}}
\newcommand{\cD}{\mathcal{D}}
\newcommand{\cE}{\mathcal{E}}
\newcommand{\cF}{\mathcal{F}}
\newcommand{\cH}{\mathcal{H}}
\newcommand{\cM}{\mathcal{M}}

\newcommand{\cL} {\mathcal{L}}

\newcommand{\Cc}{C_{\mathrm{c}}}

\input{xy}
\xyoption{all}

\newcommand{\beq}{\begin{equation}}
\newcommand{\eeq}{\end{equation}}

\newcommand{\lie}{\mathfrak}

\newcommand{\SU}{\mathrm{SU}}
\newcommand{\SL}{\mathrm{SL}}
\newcommand{\Mult}{\mathsf{Mult}}
\newcommand{\Diag}{\mathsf{Diag}}
\newcommand{\sfq}{\mathsf{q}}
\newcommand{\sfe}{\mathsf{e}}

\newcommand{\sfs}{\mathsf{s}}
\newcommand{\Ghat}{\widehat{G}}
\newcommand{\Flip}{\lie{S}}
\newcommand{\counit}{\varepsilon}
\DeclareMathOperator{\Tr}{Tr}
\newcommand{\op}{\mathrm{op}}
\newcommand{\hit}{\!\rhd\!}
\newcommand{\hitby}{\!\lhd\!}

\newcommand{\UqRsu}{\cU_q^{\R}(\fsu(2))}
\newcommand{\UqRk}{\cU_q^{\R}(\lie{k})}
\newcommand{\UqRg}{\cU_q^{\R}(\lie{g})}

\newcommand{\ip}[1]{\langle #1 \rangle}

\newcommand{\half}{{\frac12}}
\newcommand{\directint}{\int^\oplus}

\usepackage[dvipsnames]{xcolor}
\usepackage[normalem]{ulem}	
 \newcommand\blfootnote[1]{%
  \begingroup
  \renewcommand\thefootnote{}\footnote{#1}%
  \addtocounter{footnote}{-1}%
  \endgroup
}
 
\begin{document}

\Large

\centerline{\bf Quantized semisimple Lie groups}

\bigskip

\medskip
\centerline{\Large{{ R. Fioresi$,^{1,2}$ R. Yuncken$^{3}$}}}

\normalsize

\blfootnote{
These notes are based on a short series of lectures given by
R.Yuncken in Prague in September 2023. Both authors thank Charles U.
for the hospitality and the COST action CaLISTA CA21109 for support.
R.F. thanks Gnsaga-Indam, HORIZON-MSCA-2022-SE-01-01 CaLIGOLA,
MSCA-DN CaLiForNIA - 101119552,
PNRR MNESYS, PNRR National Center for HPC, Big Data and Quantum Computing.
R.Y. is supported by the French Agence Nationale de le Recherche project OpART, ANR-23-CE40-0016.
}

\bigskip

\medskip
\centerline{$^{1}$ Fabit,
  via S. Donato 15, University of Bologna, 40126 Bologna, Italy}

\medskip
\centerline{$^{2}$ Istituto Nazionale di Fisica Nucleare, Sezione di Bologna,
  40126 Bologna, Italy.}

\medskip
\centerline{{rita.fioresi@unibo.it}}

\medskip
\centerline{$^{3}$ Université de Lorraine, CNRS, IECL, F-57000 Metz, France}

\section{Introduction}

The purpose of this work
is to give a quick and basic introduction to the
quantization of semisimple Lie groups,
{specifically compact and complex semisimple Lie groups},
from the point of view of unitary representation theory.

\medskip
It is an extraordinary fact that the simple Lie algebras
over $\C$ can be completely classified.  The complete list consists of:
\begin{itemize}
\item four infinite families, denoted $A_n$, $B_n$, $C_n$ and $D_n$  ($n\in\N$), corresponding to the classical matrix groups $\mathrm{SL}(n+1,\C)$, $\mathrm{SO}(2n+1,\C)$, $\mathrm{Sp}(n,\C)$ and $\mathrm{SO}(2n,\C)$, respectively,
\item five exceptional Lie algebras, corresponding to the exceptional simple Lie groups denoted $E_6$, $E_7$, $E_8$, $F_4$ and $G_2$.
\end{itemize}

In the 1980s,  Russian mathematical physicists discovered that the
universal enveloping algebra of $\lie{sl}_2(\C)$ admits a
$q$-deformation, or quantization, as a Hopf algebra
$\cU_q(\lie{sl}_2)$, with $q\in\C^\times$ being a complex parameter \cite{Sklyanin:Uqsl2, KulRes}.
This led to an explosion of discoveries of quantum groups, most
notably the quantized enveloping algebras of Drinfeld and Jimbo  \cite{Drinfeld:quantum_groups, Drinfeld:quantum_groups2,Jimbo:Uqg},  and then
the quantized algebras of functions on matrix groups of Reshetikhin, Takhtajan and Faddeev \cite{FadResTak} and
 Woronowicz \cite{Woronowicz:SUq2,Woronowicz:pseudogroups}. For a complete account of the beginning of
such discoveries see the excellent historical notes at the ends of the chapters in \cite{KliSch}.

The moral to be taken from their discoveries is the following:
semisimple Lie groups all come  with quantum deformations, and these
quantum deformations can be observed from multiple perspectives.  It
is natural to explore these quantizations for numerous reasons,
{the most obvious being} (1)
with a view to applications in physics, as their
origin from the theory of the $R$-matrix, integrable systems and Yang
Mills theories clearly shows, and
(2) with the hope to learn more about the classical Lie groups from
which they emanate, deepening our understanding of their
representation theory.
 
For simplicity, in these notes we will concentrate almost entirely on
the rank-one simple Lie algebra
$\lie{sl}_2(\C)$
and its associated compact and complex Lie groups, $\mathrm{SU}(2)$ and
$\mathrm{SL}(2,\C)$.
Nonetheless, the results we recount can all be extended to higher-rank compact and complex semisimple groups if one is willing to invest in the structure theory of semisimple groups and their somewhat more complicated quantum analogues.  

Perhaps the biggest gap in the present state of quantized semisimple Lie groups is the question of quantizing noncompact real groups, such as $\mathrm{SL}(n,\R)$ or $\mathrm{SU}(p,q)$ from the point of view of operator algebras.  Progress is being made \cite{KoeKus:SUq11,DeCommer:SUq11, DeCDzo:SL2R}, but there are still many issues to be uncovered.

\section{{Quantized enveloping algebras}}

In this section we briefly recap a few facts on semisimple Lie algebras and
we introduce the quantization of their enveloping algebras, together with some
standard facts about their representation theory.

\subsection{Semisimple Lie Algebras}

We begin with a brief introduction to semisimple Lie algebras,
sending the reader to \cite{Varadarajan:Lie_groups}, \cite{Knapp:Lie_groups}, \cite{Helgason:book}
for the full account. We assume the reader is familiar with the notion
of Lie algebra and Lie group and their adjoint representations, denoted,
as usual, with $\ad$ and $\Ad$ respectively.  

\medskip
A complex Lie algebra {$\lie{g}$} is \textit{semisimple} if the bilinear
form (Cartan-Killing form)
$$
(X,Y):=\tr(\ad(X)\ad(Y))
$$
is non degenerate. Complex semisimple Lie algebras are direct sums of simple
ones, namely Lie algebras with no non trivial ideals
{(excluding the one-dimensional abelian Lie algebra $\C$)}. Simple Lie algebras
are completely classified and they either belong to the infinite families:
$$
A_n=\fsl_{n+1}(\C), \qquad B_n=\fso_{2n+1}(\C),
\qquad C_n=\fsp_{n}(\C), \qquad D_n=\fso_{2n}(\C)
$$
or the so called exceptional Lie algebras:
$$
G_2, \quad F_4, \quad E_6, \quad E_7, \quad E_8.
$$
As usual, $\fsl_n(\C)$ denotes the complex special linear Lie algebra,
$\fso_n(\C)$ the complex orthogonal Lie algebra and  $\fsp_{n}(\C)$
the complex symplectic Lie algebra.

\medskip
The index of the Lie algebras in the above notation represents its rank, namely the dimension
of any of its Cartan subalgebras, or CSAs, 
which are, by definition, maximal abelian
subalgebras consisting of semisimple elements. {The CSAs are all isomorphic.}
Once we
fix $\fh$ a CSA of a complex simple Lie algebra ${\fg}$,
{by diagonalizing} its action on ${\fg}$ {via} the
bracket we obtain the root space decomposition of ${\fg}$:
\beq\label{root-dec}
{\fg}=\fh \oplus \bigoplus_{\al \in \Delta} {\fg}_\al
\eeq
where ${\fg}_\al=\{X \in \fg\,|\, [H,X]=\al(H)X, \, \forall H \in \fh\}$.
The non zero $\al \in \fh^*$ appearing in (\ref{root-dec})
are called {\it roots}
and their set $\Delta$ is called the {\it root system} of $\fg$. It plays a
key role in the above mentioned classification result: two simple Lie algebras
are isomorphic if and only if they have the same root system.

\medskip
If $\fg$ is a real Lie algebra we call
$$
\fg_\C=\fg \underset{{\R}}{\otimes}\C
$$
its complexification. If $\fg_\C$ is semisimple, we say that
$\fg$ is semisimple.

\medskip
We define the {\it universal enveloping algebra} of a (real or complex)
Lie algebra $\fg$ as:
$$
\cU(\fg):=T(\fg)/I, \qquad I={\langle}[X,Y]-X\otimes Y+Y \otimes X \mid X,Y \in \fg{\rangle}
$$
where $T(\fg)$ denotes the tensor algebra over $\fg$ and $I$ is
the two sided ideal {generated by all elements of the given form}.
Thus
$\cU(\fg)$ is an associative algebra and there is a one to one correspondence
between the representations of $\fg$ and the representations of its
universal enveloping algebra.

\medskip
Let us look at an interesting example.

\begin{example}
\label{ex:su2}
Consider the real Lie algebra $\fso(3)$, that is the Lie
algebra of skew symmetric matrices:
$$
\fso(3)=\{A \in M_3(\R)\,|\, A=-A^t\}
$$
where $M_3(\R)$ denotes the $3 \times 3$ real matrices. It is generated
both as vector space and as a  Lie algebra by the matrices
$X=E_{23}-E_{32}$, $Y=-E_{13}+E_{31}$, $Z=-E_{12}+E_{21}$. We leave
as a simple exercise to the reader to verify that
$$
[X,Y]=Z, \qquad [Y,Z]=X, \qquad [Z,X]=Y
$$
This real Lie algebra is isomorphic to the special unitary
Lie algebra of skew hermitian matrices:
$$
\fsu(2)=\{A \in M_2(\C)\,|\, A=-A^*\}
$$
where $A^*=\overline{A}^t$.
{For instance, a Lie algebra isomorphism is obtained by
$$
  X \mapsto \frac12 \begin{pmatrix} 0& i \\ i&0 \end{pmatrix}, \quad
  Y \mapsto \frac12 \begin{pmatrix} 0& -1 \\ 1&0 \end{pmatrix}, \quad
  Z \mapsto \frac12 \begin{pmatrix} i& 0 \\ 0&-i \end{pmatrix}, \quad  
$$
}

As one can readily check:
$$
\fsl_2(\C)=\fsu(2) \underset{{\R}}{\otimes} \C
$$
where $\fsl_2(\C)$ is the Lie algebra of $2\times 2$-matrices with zero trace. 
For this complexified Lie algebra, rather than use the basis $X,Y,Z$
above, it is better to use the basis
$$
H=\begin{pmatrix} 1 & 0 \\ 0 & -1 \end{pmatrix}, \quad
E=\begin{pmatrix} 0 & 1 \\ 0 & 0 \end{pmatrix}, \quad
F=\begin{pmatrix} 0 & 0 \\ 1 & 0 \end{pmatrix}.
$$
The elements $H$, $E$ and $F$ have brackets:
\begin{equation}
\label{eq:sl2-relations}
[H,E]=2E, \quad [H,F]=-2F, \quad [E,F]=H.
\end{equation}
The universal enveloping algebra is {the complex associative}
algebra generated by the elements $H$, $E$, $F$ subject to the
{relations \eqref{eq:sl2-relations}, where they are interpreted as commutator brackets,
  \emph{i.e.},}
\begin{equation}
\label{eq:Usl2}
\cU(\fsl_2)=\langle H,E,F \,|\,
HE-EH=2E, \, HF-FH=-2F, \, EF-FE=H\rangle
\end{equation}
    {As is common, we are writing $\cU(\fsl_2)$ to denote $\cU(\fsl_2(\C))$,
      with the base field assumed to be $\C$ unless otherwise specified.}
\end{example}

\subsection{Quantized enveloping algebras}

{Around 1980, physicists studying quantum scattering theory, notably Kulish, Reshetikin \cite{KulRes} and Sklyanin \cite{Sklyanin:Uqsl2} observed that a variant of the enveloping algebra $\cU(\fsl_2)$ appears when finding solutions to matrix equations which arise in the quantum inverse scattering method.  Informally, this quantized enveloping algebra can be presented as follows:}  
\begin{equation}
\label{eq:Uqsl2}
\cU_q(\fsl_2)=\langle H,E,F \,|\,
[H,E]=2E, \, [H,F]=-2F, \, [E,F]=[H]_q\rangle
\end{equation}
The only change from the classical enveloping algebra $\cU(\fsl_2)$ in Equation \eqref{eq:Usl2} is the replacement of $H$ by the expression $[H]_q$, which remains to be explained.  Before we get to this, we need to introduce $q$-numbers.

\begin{observation}
Let $q$ be a real positive number and let $a \in \C$. We define:
$$
[a]_q=\frac{q^a-q^{-a}}{q-q^{-1}}
$$

In the special case where $a$ is a natural number, one
can readily see that {this is a Laurent polynomial in $q$}:
$$
[a]_q=q^{a-1}+q^{a-3}+ \dots q^{-a+1}
$$
As $q$ tends to 1 we have $[a]_q \lra a$.

These $q$-numbers have some properties resembling the integers, for instance:
$$[a-b]_q[a+b]_q=[a]_q^2-[b]^2_q.$$
{On the other hand, of course, $[a]_q + [b]_q \neq [a+b]_q$.}
\end{observation}

{With this definition in hand, we wish to make sense of the formal expression
\[
[H]_q=\frac{q^H-q^{-H}}{q-q^{-1}},
\]
where $H$ is an element of an algebra.

\begin{observation}\label{H_vs_qH} 
There are three different ways to achieve this, according to one's preferences:

\begin{enumerate}

\item By means of a formal power series in $q=e^h$, where $h$ is an
indeterminate. In this point of view, we define $\cU_q(\fsl_2)$
to be an algebra over the
ring of formal power series $\C[[h]]$, with generators $E,F,H$, and we write
$$
[H]_q=\frac{e^{hH}-e^{-hH}}{e^h-e^{-h}},
$$
where $e^{hH}$ denotes the usual exponential power series:
$$
1+hH+\frac{1}{2!}h^2H^2 + \cdots
$$

This point of view is favored by algebraists and those working in formal deformation theory.
We will not use it.

\item By functional calculus, after representing the algebra 
generated by $E$, $F$ and $H$ in a suitable Hilbert space. {Once $H$ is realized as a self-adjoint operator on a Hilbert space, the expression $[H]_q=\frac{q^H-q^{-H}}{q-q^{-1}}$, makes perfect sense.

This approach is less satisfactory from a formal point of view, since it requires an implicit understanding of the representation theory of the algebra $\cU_q(\fsl_2)$ prior to giving its rigorous definition.

\item By 
replacing the generator $H$ by} a new generator $K=q^H$ {and its inverse $K^{-1} = q^{-H}$.  
In this solution, we must alter the list of relations in Equation
\eqref{eq:Uqsl2}
to eliminate any reference to $H$.  We thus define}
 $\cU_q(\fsl_2)$
as the associative algebra generated by $E$, $F$ and $K^{\pm 1}$ and
subject to the relations:
\begin{equation}
\label{eq:EFK-relations}
KEK^{-1}=q^{-1}E, \quad KFK^{-1}=q^{-2}F, \quad
[E,F]=\frac{K-K^{-1}}{q-q^{-1}}
\end{equation}
\end{enumerate}
\end{observation}

The following observation, whose proof we leave to the reader,
explains why the first two relations
in \eqref{eq:EFK-relations} are the appropriate
replacements of the first two relations in \eqref{eq:Uqsl2}.

\begin{observation}
Fix $q>0$.  Let $H$ and $X$ be operators on a finite dimensional Hilbert space, with $H$ self-adjoint, and let $\lambda\in\R$.  The following are equivalent:
\begin{enumerate}
 \item
 $[H,X] = \lambda X$
 
 \item
  $q^{H} X q^{-H} = q^\lambda X$ for any $q>0$.  Here $q^H$ is understood as an operator defined via functional calculus, using $q^H=e^{hH}$ with $h=\ln(q)$.
\end{enumerate}
\end{observation}

\subsection{Representations of $\cU_q(\fg)$}\label{sec:highest weight modules}
\label{sec:Vm}

We are interested in the finite dimensional representations
of the quantized universal enveloping algebra $\cU_q(\fg)$,
$\fg=\fsl_2(\C)$ introduced in our previous section.

Since $H$ is semisimple (i.e. diagonalizable), it acts diagonally
on $\cU(\fg)$ with real spectrum. Notice that this fact is immediately
generalized to the case of a CSA $\fh$ of a complex simple
Lie algebra, since $\fh$ consists of commuting semisimple elements.

We now focus on the quantization {of} $\cU(\fg)$.

\begin{definition}
We say that a finite dimensional complex representation of $\cU_q(\fg)$ is 
\textit{integrable} if $q^H$ acts diagonally with spectrum consisting
of positive real numbers.
\end{definition}

We have the following classification result, mimicking a
corresponding classical result, i.e. in the case $q=1$, {which can be} found
for example in \cite{Varadarajan:Lie_groups} Ch. 4.

\begin{theorem} 
\label{thm:Vm}
Every finite dimensional integrable representation of
$\cU_q(\fg)$ is isomorphic, for a suitable $m \in \half\N$ to
$$
V(m)=\mathrm{span}\{v_m,v_{m-1}, \dots v_{-m}\}
$$
with
$$
\begin{array}{c}
Hv_{\mu}=2\mu v_\mu, \qquad
Fv_{\mu}=v_{\mu-1},\qquad Ev_\mu=[m-\mu]_q[m+\mu+1]_qv_{\mu+1}
\end{array}
$$
{for all $\mu$, with the convention $v_{\mu}=0$ if $\mu\notin\{m,m-1,\ldots,-m\}$.}
\end{theorem}

We call $m, m-1 \dots$ the \textit{weights} of the representation
and $v_m, v_{m-1}, \dots$ are \textit{weight vectors}.

\begin{remark}
Note that we are using the convention that our weights are representated \emph{half} integers\footnote{In \cite{Varadarajan:Lie_groups}, what we call a weight vector of weight $k$
is called a weight vector of weight $2k$. All formulae in there have to be interpreted accordingly.}, $\mu\in \half\N$.
This is reflected in the factor of $2\mu$ in the  action of $H$ given in the theorem, and it will have consequences in certain formulas to follow.  There is an equally common convention where it is the whole integer value $2\mu\in\N$ which is called the weight.  
\end{remark}

\begin{proof}
Existence. We define $V(m)$
as above and then verify the relations in \eqref{eq:Uqsl2}, or equivalently \eqref{eq:EFK-relations}.
We leave
this check as an exercise.

Uniqueness. Assume we have a finite dimensional irreducible
complex representation
$V$ of $\cU_q(\fg)$. Let $m$ be the {highest} weight, i.e.\  {half} of the eigenvalue  {of $H$}
with largest real part, and let $v_m$ be one of its eigenvectors.

We have the following important fact. If $v \in V$ has weight $\lambda$
then
$$
HEv=[H,E]v+E(Hv)=2Ev+E(2\lambda v)=2(\lambda +1)E.
$$
Hence we have that $Ev$ is a weight vector of weight $\lambda+1$.
Similarly $Fv$ is a weight vector of weight $\lambda-1$.
The consequence of this observation is that
$$
\mathrm{span}\{v_m,v_{m-1}, \dots v_{-m}\}
$$
is invariant by the action of $H$, $E$ and $F$.  By
the irreducibility hypothesis it coincides with $V$.

As an exercise one can prove that
 {
\beq \label{dimV}
[E,F^{k+1}]=[k+1]_q[H+k]_qF^{k}.
\eeq
}
Since $V$ is finite dimensional we must have $F^{k+1}v_m=0$ for some $k\in\N$,  {and taking the smallest such $k$ we have $F^kv_m\neq 0$.  Since  $Ev_m=0$, we then have $[E,F^{k+1}]v_m=0$ and hence $[H+k]_qF^kv_m=0$.  Thus $F^kv_m$ is in the kernel of $H+k$, and since $F^kv_m$ has weight $m-k$, this means $2(m-k)+k=0$, or $k=2m$.  Therefore $m\in\half\N$.  Putting $v_{m-k} = F^kv_m$, we obtain a basis for $V(m)$ with the stated actions of $H$ and $F$.  The action of $E$ on the $v_\mu$ can be calculated from Equation \eqref{dimV}.
}
\end{proof}

We end this section with {some remarks on} Verma modules. For more
details see \cite{Dixmier:Enveloping_algebras}.

\begin{observation}
\label{obs:Verma}
 {Let now $m\in\C$.}  
Let us define  {a countably infinite dimensional} $\cU_q(\fg)$ representation:
$$
M(m)=\mathrm{span}\{v_m,v_{m-1}, v_{m-2} \dots \}
$$
subject to the action of $H$, $E$, $F$   {with exactly the same formulas as in Theorem \ref{thm:Vm}}.
{This module is called the \emph{Verma module} of highest weight $M(m)$.  Notice that $M(m)$ is not finite dimensional, nor necessarily irreducible}.

One can show that if $m \not\in\half\N+i\pi h^{-1}\Z$, then
$M(m)$ is irreducible.  {On the other hand, if $m \in \half\N+i\pi h^{-1}\Z$, then} $M(m)$ is reducible  {and} there is an exact sequence:
 {
\begin{equation}
\label{eq:BGG1}
0 \lra M(-m-1) \hookrightarrow M(m) \lra V(m) \lra 0
\end{equation}
}
{This resolution of the irreducible finite dimensional modules $V(m)$ by Verma modules plays an important role in many parts of representation theory.}
\end{observation}

\section{Hopf algebras and quantum groups}

In this section we introduce our main object of study: Hopf algebras.  {Hopf algebras are the central algebraic objects in all approaches to quantum groups.  As we will see, quantized enveloping algebras, algebras of functions on quantum groups, and convolution algebras of quantum groups will all be modelled by Hopf algebras.}

For more details on definitions and techniques see
\cite{Montgomery:Hopf_lectures,KliSch}.

\subsection{{Motivation: Some rough remarks on Pontryagin duality}}
\label{sec:Pontryagin}

{In this subsection, which is not intended to be {complete},
we make some general remarks on the duality
between a given group, or more precisely its group algebra,
and the functions on such a group. We want to give some hint to
the so called Pontryagin duality, a deep result, which is a key tool for our
understanding of quantum geometry.}

Consider the example of the Lie group $\R$.  
Recall that the Fourier transform gives an isomorphism between the algebra of functions on $\R$ with pointwise multiplication, and the algebra of functions on $\R$ with convolution.  To make this precise requires some careful analysis of the spaces of functions.  Several specific algebra isomorphisms are possible, for instance:
\begin{itemize}
\item the Schwartz space $\mathcal{S}(\R)$ of rapidly decaying functions with pointwise multiplication, and  $\mathcal{S}(\R)$ with convolution;
\item the algebra of continuous functions vanishing at infinity $C_0(\R)$ with pointwise multiplication, and the group $C^*$-algebra $C^*(\R)$ with convolution;
\item the Paley-Weiner-Schwartz isomorphism between functions on $\R$ which extend analytically to $\C$ 
with certain growth conditions under pointwise multiplication, and compactly supported distributions on $\R$ with convolution;
\item the Fourier algebra $\mathrm{A}(\R)$ with pointwise multiplication, and the algebra $L^1(\R)$ with convolution.
\end{itemize}

Similar isomorphisms can be obtained between functions on the circle $\mathbb{T}$ with pointwise multiplication, and functions on $\Z$ with convolution, or vice versa, by using Fourier series instead of Fourier transforms.  

These are examples of Pontryagin duality.  We will not make explicit use of Pontryagin duality in these notes, but it can be useful to motivate the constructions to follow.  {Let us develop this philosophy without going into technical details.}

\medskip

Let $G$ be a group.  For the moment, we may consider a finite group, a
locally compact topological group, a Lie group, or a number of other
types of group.  

{As indicated above,} a group admits two general
types of algebras.
The first general type of algebra is ``functions with pointwise multiplication''.
Some examples of this type:
\begin{itemize}
\item $\C(G)$,  functions on a finite group,
\item $C_0(G)$, continuous functions vanishing at infinity on a locally compact group,
\item $\cO(G)$, regular functions on a linear algebraic group,
\item $\Cc^\infty(G)$, smooth compactly supported functions on a Lie group,
\item $L^\infty(G)$, essentially bounded functions on a measured group.
\end{itemize}
The second type is ``distributions with convolution'':
\begin{itemize}
\item $\C[G]$, the group ring of a finite group,
\item $L^1(G)$, integrable functions on a locally compact topological group (which identify as distributions upon multiplying with Haar measure),
\item $\Cc^\infty(G)$, smooth compactly supported functions on a Lie group with convolution (likewise),
\item $\cE'(G)$, compactly supported distributions on a Lie group,
\item $C^*(G)$, the group $C^*$-algebra of a locally compact group.
\item $\mathrm{vN}(G)$, the von Neumann algebra of a locally compact group,
\item $\cU(\lie{g})$, the universal enveloping algebra of a Lie algebra (which are distributions supported at the identity).
\end{itemize}
  Remaining vague about the analytical details, let us denote any algebra of the first kind by $\cA(G)$ and any of the second kind as $\cD(G)$.  We will be more concrete in the next section.

The reason why these two types of algebras exist is a consequence of the functoriality of functions and distributions with respect to the fundamental operations of a topological group.  Recall that {$\C$-valued} functions are contravariant objects.  For instance, if $\phi:M\to N$ is a smooth map between manifolds, then we can pull back {smooth $\C$-valued} functions from $N$ to $M$:
\[
  \phi^*:C^\infty(N) \to C^\infty(M); \qquad f \mapsto \phi^*f=f\circ\phi.
\]
On the other hand, distributions are covariant objects, since they are dual to smooth functions. {Explicitly, if we write $(u,f)$ for the pairing between a distribution $u\in\cE'(M)$ and a smooth function $f\in C^\infty(M)$, then we have a pushforward map defined by}
\[
  \phi_*:\mathcal{E}'(M) \to \mathcal{E}'(N); \qquad (\phi_*u,f) := (u,\phi^*f)
\]
for $u\in\mathcal{E}'(M)$, $f\in C^\infty(N)$.

{Now consider the case where $M$ is a group.}
The pointwise products and convolution products {of the various algebras} listed above all arise by applying the functorial operations to two basic structural maps on a group:
\begin{align*}
  &\Mult : G\times G \to G; \quad  (g,h) \mapsto  gh &&\text{(the group law)} \\
  &\Diag : G \to G\times G; \quad g \mapsto (g,g) &&\text{(the diagonal embedding)}.  
\end{align*}
The diagonal embedding is a greatly underestimated map.  We shall see its crucial role in Pontrjagin duality shortly.

These basic structural maps induce products by push-forward and pull-back:
\begin{align*}
  & \Mult_* :  \cD(G) \otimes \cD(G) \cong \cD(G\times G) \to \cD(G) && \text{(convolution)} \\
  & \Diag^*:   \cA(G) \otimes \cA(G) \cong \cA(G \times G) \to \cA(G) && \text{(pointwise multiplication)}.
\end{align*}
The isomorphisms here will require using some completed tensor product
$\otimes$ which will depend heavily upon the categories of
``functions'' and ``distributions'' in which one is working.
Since we are deliberately avoiding the analytical details in this section, we {will not} address this issue here.

This point of view reveals that the algebras $\cA(G)$ and $\cD(G)$ should also have coproducts, \emph{i.e.} maps in the reverse direction:
\begin{align*}
  & \Diag_*:   \cA(G) \to  \cA(G \times G) \cong \cA(G) \otimes \cA(G) , \\
  & \Mult_* :  \cD(G) \to \cD(G\times G) \cong \cD(G) \otimes \cD(G) . \\
\end{align*}
Thus $\cA(G)$ and $\cD(G)$ will be \emph{bialgebras}, and more specifically \emph{Hopf algebras} (see below for definitions).  Again, for infinite groups, the technicalities of the tensor product may require us to add some qualifying adjectives---\emph{multiplier Hopf algebras}, for instance.  

The principle of Pontrjagin duality is that one should be able to exchange
the roles of the two algebras $\cA(G)$ and $\cD(G)$.
For instance, for the group $G=\R$ we discussed various instances of algebra isomorphisms of the general form $\cA(\R) \cong \cD(\R)$ for different examples of pointwise algebras $\cA(\R)$ and convolution algebras $\cD(\R)$.  We say that $\R$ is its own Pontryagin dual, and write $\hat{\R} = \R$.  Similarly, there are various isomorphisms of the form $\cA(\T) \cong \cD(\Z)$, giving the Pontryagin duality $\hat{\T} = \Z$.

\medskip
It is natural to ask whether Pontryagin duality applies more
generally. To begin with, we could
formally define an algebra $\cA(\Ghat) := \cD(G)$ without giving any
concrete meaning to the symbol $\widehat{G}$.  If we are lucky, we will be able to realize this algebra $\cA(\Ghat)$ as the algebra of functions on some concrete topological group $\Ghat$.

However, this strategy relies on the categories of algebras $\cD(G)$
and $\cA(G)$ being the same,
which is not the case.
{The fundamental obstruction is that} the pointwise algebras $\cA(G)$ are always abelian, whereas
the convolution algebras
$\cD(G)$ are abelian {only} if $G$ is.  Let us spell this out.

The properties of the algebras $\cA(G)$ and $\cD(G)$ are inherited from properties of $\Diag$ and $\Mult$, respectively.  For instance, the associativity of $\cD(G)$ comes from the associativity of $\Mult$,
\begin{align*}
  \Mult \circ (\Mult \times \id)= \Mult \circ (\id\times\Mult) : \quad (g,h,k) \mapsto ghk,
\end{align*}
while the associativity of $\cA(G)$ comes from the co-associativity of $\Diag$,
\[
  (\Diag \times \id) \circ \Diag = (\id \times \Diag) \circ \Diag : \quad g \mapsto (g,g,g).
\] 
Similarly, $\cA(G)$ is always abelian, since $\Diag$ is always co-commutative in the sense that if $\Flip:(g,h) \mapsto (h,g)$ denotes the flip map on $G\times G$, we have
\[
  \Flip\circ\Diag = \Diag,
\]
On the other hand, $\cD(G)$ is abelian only if $\Mult$ is abelian, that is, only if
\[
  \Mult\circ\Flip = \Mult.
\]

If we hope to have a coherent theory of Pontrjagin duality, we must make one of the following choices:
\begin{enumerate}
\item Restrict attention to abelian groups for which $\Mult$ is commutative.
\item Define ``non co-abelian groups'' for which $\Diag$ is not cocommutative.
\end{enumerate}
The second strategy, appropriately interpreted, leads to the definition of a quantum group.  The {specific technical} details will depend upon the category of classical groups one is trying to mimic,
{but will in every case require the notion of Hopf algebra.}

\subsection{Hopf algebras}

Let us make some of the last section more concrete.  Here is the quick
definition of
a Hopf algebra, which is spelled out in more detail in the Appendix \ref{Hopf-app}.

\begin{definition}
  A Hopf algebra over $\C$ is a complex vector space $\cA$ equipped with five maps:
\begin{itemize}
\item An associative product ${\mu} : \cA \otimes \cA \to \cA$,
\item A coassociative coproduct $\Delta : \cA \to \cA \otimes \cA$,
\item A unit $1\in\cA$, realized also as a unit map  ${i}:\C\to \cA$; $z\mapsto z1$,
\item A counit map $\counit:\cA \to \C$, satisfying $(\id\otimes\counit)\circ\Delta = \id = (\counit\otimes\id)\circ\Delta$,
\item An invertible antipode $S:\cA \to \cA$, verifying the antipode identity
\[
  \mu\circ (S\otimes\id) \circ \Delta = {i} \circ \counit =  \mu\circ (\id\otimes S) \circ \Delta.
\]
\end{itemize}
The coproduct and counit should be unital algebra homomorphisms, or equivalently the product and unit are counital coalgebra homomorphisms.  The antipode is automatically an algebra and coalgebra antihomomorphism.
\end{definition}

We shall use Sweedler notation for the coproduct:
\begin{equation}
\label{eq:Sweedler}
  \Delta (X) = X_{(1)} \otimes X_{(2)} ,
\end{equation}
where the right-hand side is a formal notation to designate a finite sum of elementary tensor products $\sum_i X_{(1),i} \otimes X_{(2),i}\in \cA \otimes \cA$.

The following examples show that finite dimensional Hopf algebras
{can be used to give rigorous definitions of algebras of both types from the previous section,
$\cA(G)$ and $\cD(G)$, in the case where $G$ is a finite group.}

Let $G$ be a finite group, with the following five structural maps
\begin{align*}
\Mult:&\ G\times G \to G,&&\text{group law }\\
\Diag:&\ G \to G\times G, &&\text{diagonal embedding} \\
\sfe: &\ \{1\} \to G, &&\text{inclusion of the unit} \\
\sfq:&\ G \to \{1\}, &&\text{quotient to the unit} \\
\sfs:&\ G\to G, &&\text{inverse}.
\end{align*}

\begin{example}
\label{ex:finite-A}
The algebra  $C(G)$ of complex-valued functions on $G$ is a Hopf algebra with structural maps defined by pull-back: 
\[
(\mu,\Delta,{i},\counit,S) = (\Diag^*, \Mult^*, \sfq^*, \sfe^*, \sfs^*).
\]
This is the Hopf algebra of ``functions with pointwise multiplication''.
\end{example}

\begin{example}
\label{ex:finite-D}
The dual space $C(G)^*$ is a Hopf algebra with structural maps defined by push-forward:
\[
 (\mu,\Delta,{i},\counit,S) = (\Mult_*, \Diag_*, \sfe_*, \sfq_*, \sfs_*).
\]
Morally, this is the Hopf algebra of ``distributions with convolution''.
\end{example}

Example \ref{ex:finite-A} is commutative and Example \ref{ex:finite-D} is cocommutative.  There are examples of finite dimensional Hopf algebras which are neither commutative nor cocommutative.  But we will be interested in infinite dimensional examples.

\subsection{Quantized enveloping algebras as quantum groups}
In this section we present an example of a quantized enveloping algebra {as a Hopf algebra}.  {We begin with the classical case.}

\begin{example}
\label{ex:Ug}
Let $\lie{g}$ be a complex Lie algebra.  The universal enveloping algebra $\cU(\lie{g})$ is a Hopf algebra where the coproduct, counit and antipode are defined on generators $X\in\lie{g}$ by
\[
  \Delta(X) = X\otimes1 + 1\otimes X, \quad \counit(X)=0, \quad S(X) = -X.
\]

It is an exercise to show that these are precisely the maps $\Diag_*$, $\sfq_*$ and $\sfs_*$, if the elements of $\cU(\lie{g})$ are viewed as distributions supported at the identity.  
\end{example}

After Kulish and Reshetikhin \cite{KulRes} discovered the quantized enveloping algebras $\cU_q(\lie{sl}_2)$, Sklyanin \cite{Sklyanin:Uqsl2} observed that it is {also} a Hopf algebra.
Thus quantum groups were born.

For the purpose of these notes, we call \textit{quantum group} a Hopf algebra
which is related, typically by a deformation parameter $q$ as above, to either
the universal enveloping algebra, or the function algebra (Sec. \ref{fn-sec}), {or the convolution algebra (Definition \ref{def:DKq})}
of an algebraic group.

\begin{theorem}[Sklyanin]
The quantized enveloping algebra $\cU_q(\lie{sl}_2)$ with generators $E,F,K^{\pm1}$ and relations \eqref{eq:EFK-relations} is a Hopf algebra when equipped with the $q$-deformed maps
\[
  \Delta(E) = E \otimes q^H + 1 \otimes E, \quad \Delta(F) = F\otimes 1 + q^{-H} \otimes F, \quad \Delta(q^H) = q^H\otimes q^H,
\]
\[
  \counit(E) = \counit(F) = 0, \quad \counit (q^H) = 1,
\]
\[
  S(E) = -Eq^{-H}, \quad S(F) = -q^HF, \quad S(q^H) = q^{-H}.
\]
\end{theorem}

\begin{observation}
The formulas for $\Delta(q^H)$, $\counit(q^H)$ and $S(q^H)$ follow from the classical formulas
\[
  \Delta(H) = H\otimes 1 + 1\otimes H, \quad \counit(H)=0, \quad S(H) = -H,
\]
via functional calculus, in the same way as in Observation \ref{H_vs_qH}
--- exercise for the reader.  In this sense, it is reasonable to say that the generator $H$ is undeformed in $\cU_q(\fsl_2)$.  On the other hand, the formulas for $E$ and $F$ are $q$-deformed, and putting $q=1$ recovers the classical formulae.  
\end{observation}

\subsection{Real structures {and compact quantum groups}}
\label{sec:real_structures}

A real structure on a quantum group is encoded by a $*$-structure.  To understand why this is the case, let us begin with the classical case.

Let $K$ be a \emph{real} Lie group, with real Lie algebra $\lie{k}$.  For simplicity we take $K$ to be compact {here}.  

  A unitary representation of $K$ is a continuous homomorphism $\pi:K\to U(V)$ from $K$ to the unitary operators on a complex Hilbert space $V$, which we will assume is finite dimensional.  By differentiating $\pi$ at the identity, we obtain a Lie algebra homomorphism from $\lie{k}$ to the Lie algebra $\lie{u}(V)$ of skew-Hermitian operators on $V$.  This is again denoted by $\pi$.  
  
Note that the Lie algebra representation $\pi:\lie{k} \to \lie{u}(\lie{k})$ is a real linear map.  As usual, it is far preferable to work with complex linear maps, since spectral theory in complex vector spaces is better behaved.  We are therefore led to define the complex linear extension of $\pi$.  Specifically, if we write
\[
  \lie{g} = \lie{k}_\C = \lie{k} \otimes_\R \C
\]
for the complexification of $\lie{k}$, then we can extend $\pi$ by complex linearity to a complex Lie algebra homomorphism
\[
  \pi : \lie{g} \to \lie{u}(V)_\C = \End(V).
\]
As is the habit in this subject, we continue to denote this complexified representation by the same letter $\pi$.  
Finally, we can extend this complex Lie algebra morphism {by universality} to a complex algebra homomorphism from the universal enveloping algebra,
\begin{equation}
\label{eq:pi}
  \pi: \cU(\lie{g}) \to \End(V).
\end{equation}
It is this notion of representation that is most easily carried over to the quantized enveloping algebras $\cU_q(\lie{g})$.

Note, though, that a homomorphism of this kind could also be obtained from a different starting point.  Let $G$ denote the complex Lie group whose Lie algebra is $\lie{g}$.  Suppose that $\sigma:G\to\End(V)$ is a \emph{holomorphic} (not unitary) representation of $G$ on $V$.  We can pass to the complex derivative of $\sigma$ at the identity and then extend to a complex algebra representation
\begin{equation}
\label{eq:sigma}
  \sigma:\cU(\lie{g}) \to \End(V).
\end{equation}

This suggests the question: how can we distinguish a holomorphic representation {\eqref{eq:sigma}} of $\lie{g}$ from a unitary representation {\eqref{eq:pi}} of $\lie{k}$, once we have passed to the associated representation of $\cU(\lie{g})$?  The answer is that unitary representations of $\lie{k}$ have one further piece of structure.   Namely, since we can distinguish the real linear subspace $\lie{k}$ inside $\lie{g}$, we can define a $*$-involution on $\lie{g}$:
\begin{equation}
\label{eq:star}
  * : X \mapsto 
  \begin{cases}
   -X &\text{if } X\in\lie{k}, \\
   X & \text{if } X \in i\lie{k}.
  \end{cases}
\end{equation}
This definition is motivated by the fact that under any unitary representation of $K$, elements of $\lie{k}$ will act as skew-Hermitian operators.

The involution \eqref{eq:star} is a complex antilinear anti-homomorphism of $\lie{g}$, meaning
\begin{align*}
  (\lambda X)^* &= \overline{\lambda}X^*, \\
   [X,Y]^* &= [Y^*,X^*]
\end{align*}
for all $\lambda\in\C$ and $X,Y\in\lie{g}$.  It extends to a complex antilinear algebra antihomomorphism $*$ of $\cU(\lie{g})$.  That is, $\cU(\lie{g})$ becomes a $*$-algebra.  

When $\cU(\lie{g})$ is equipped with the involution $*$ of \eqref{eq:star}, we will denote it instead by $\cU^\R(\lie{k})$, to signify the fact that we are regarding {it}
as the complexification of the enveloping algebra of a \emph{real} Lie group $K$, and not as the enveloping algebra of the complex Lie group $G$.  

\begin{example}
Consider the example of $\lie{k}=\lie{su}(2)$.  As discussed in Example \ref{ex:su2}, the complexification of $\lie{k}$ is $\lie{g}=\lie{sl_2}(\C)$.  The standard basis elements $E,F,H$ of $\lie{sl}_2(\C)$ do not belong to the real Lie subalgebra $\lie{su}(2)$.  If we write them in terms of elements of $\lie{su}(2)$, we have
\begin{align*}
  E = \begin{pmatrix} 0&1\\0&0 \end{pmatrix} &= -Y -iX, &
  F = \begin{pmatrix} 0&0\\1&0 \end{pmatrix} &= Y-iX, \\
  H = \begin{pmatrix} 1&0\\0&-1 \end{pmatrix} &= -2iZ,
\end{align*}
where 
\[
  X= \half \begin{pmatrix} 0&i\\i&0 \end{pmatrix} , \qquad
  Y= \half \begin{pmatrix} 0&-1\\1&0 \end{pmatrix} , \qquad 
  Z= \half \begin{pmatrix} i&0\\0&-i \end{pmatrix} , \qquad  
\]
form an $\R$-basis for $\lie{su}(2)$.  Applying the $*$-operation of \eqref{eq:star}, we get
\begin{align*}
  E^*& = Y -iX = F, &
  F^* &= -Y-iX = E,  &
  H^* &= -2iZ = H
\end{align*}
\end{example}

This real structure, which distinguishes $\cU^\R(\lie{su}(2))$ from $\cU(\lie{sl}_2)$, can be generalized to define the quantization $\cU_q^\R(\lie{su}(2))$, {as we will now explain}.  We start with the definition of a real structure on a {general} Hopf algebra.

\begin{definition}
A Hopf $*$-algebra is a Hopf algebra equipped with a conjugate-linear involution $*$ which is an algebra anti-automorphism and coalgebra automorphism.  
\end{definition}

\begin{proposition}
\label{prop:UqRk}
The quantized enveloping algebra $\cU_q(\fsl_2)$ is a Hopf $*$-algebra
when equipped with the involution defined by
  \[
    E^* = q^HF, \qquad F^* = Eq^{-H}, \qquad H^*=H.
  \]
\end{proposition}

\begin{remark}
 Note that, when considering the generator $q^H$, the relation $H^*=H$ is equivalent to $(q^H)^*=q^H$.  Recall that we are only considering positive real values of $q$.
\end{remark}

We will write $\UqRsu$ for  $\cU_q(\fsl_2)$ equipped with this real structure.
 
With this in hand, a \emph{unitary representation} of the quantum group $K_q$ should be understood to mean an integrable $*$-representation of the $*$-Hopf algebra $\UqRk$.
The following proposition shows that, as in the classical case, {all of the irreducible $\cU_q(\lie{g})$-modules $V(m)$ of Theorem \ref{thm:Vm}} can be realized as unitary $K_q$-representations.
 
 \begin{proposition}
   Each of the irreducible representations $V(m)$ of $\cU_q(\lie{g})$, with $m\in\half\N$, admits an inner product such that it becomes a $*$-representation of $\UqRk$.
 \end{proposition}

\subsection{Operations on Hopf algebra representations}

One of the significant advantages of Hopf algebras is that we can define tensor products of their representations.  There is no reasonable way to define the tensor product $V\otimes W$ of two representations of a general associative algebra $\cA$.  However, if $\cA$ is a Hopf algebra, we can use the coproduct $\Delta:\cA \to \cA \otimes \cA$ to define
\begin{equation}
\label{eq:tensor_repn}
  X  (v\otimes w) := X_{(1)}v\otimes X_{(2)}w, \qquad (X\in \cA,\ v\in V, \ w\in W).
\end{equation}
We are again using Sweedler notation, see Equation \eqref{eq:Sweedler}.

Moreover, the counit allows us to define a \emph{trivial representation} of $\cA$ on $\C$:
\[
  X z := \epsilon(X)z, \qquad (X\in \cA, \ z\in\C).
\]
It is a simple exercise using the axioms of a Hopf algebra
to show that $\C \otimes V \cong V \cong V \otimes\C$ as representations.

Finally, the antipode lets us define a contragredient representation on $V^*$, by
\[
  X \eta = S(X)^\mathsf{t} \eta := \eta\circ S(X), \qquad (X\in\cA, \eta\in V^*),
\]
where the superscript $\mathsf{t}$ denotes transpose.

\section{Quantized algebras of functions}\label{fn-sec}

In this section we introduce the quantized algebras of functions
of a complex semisimple group $G$ and its maximal compact subgroups $K$. We elucidate
the theory as always for $G=\mathrm{SL}(2,\C)$ and $K=\SU(2)$, though it can
be formulated in full generality (see \cite{KliSch}).

\subsection{The algebra of polynomial functions on $G_q$}
{Given a vector field on a manifold $M$, or more generally a scalar-valued differential operator on $M$, we can obtain a linear functional on $C^\infty(M)$ by {applying it and} evaluating at a point.}
When the manifold is a connected algebraic group $G$,
we can look at the Lie algebra $\fg$ of left invariant vector fields on $G$. We then 
have a duality {of Hopf algebras} between the left invariant
differential operators, which are identified with the universal enveloping algebra $\cU(\fg)$,
and the algebraic functions on $G$, $\cO(G)$:
\begin{align}
\label{diff-ac}
\cU(\fg) \times \cO(G) &\lra \C, \\
 X, f &\mapsto {(} X, f {)}
:=\left.\frac{d}{dt}\right|_0 f(\exp(tX)) &&\text{ for $X\in\lie{g}$.}
\nonumber
\end{align}
The term {\sl duality {of Hopf algebras}} refers to the fact that
{the product of $\cU(\fg)$ and coproduct of $\cO(G)$ are
dual maps with respect to this pairing, and vice-versa},
see {Observation}
\ref{ex:pairing} as well as
\cite[Ch.1]{Kassel} and our Appendix \ref{Hopf-app} for more details.

\begin{remark}
Abstractly, this duality comes from the fact that the coproduct on $\cU(\fg)$ and product on $\cO(G)$ are given by $\Diag_*$ and $\Diag^*$, respectively, according to the general philosophy of Section \ref{sec:Pontryagin}.  Similarly, the product on $\cU(\fg)$ and coproduct on $\cO(G)$ are given by $\Mult_*$ and $\Mult^*$, respectively.
\end{remark}

We thereby obtain an embedding $\cO(G) \hookrightarrow \cU(\lie{g})^*$.
This embedding is injective by the fact that algebraic functions $\cO(G)$ are determined
by their jets at the identity.

\medskip
We now turn to the quantum setting and realize, with the same philosophy,
the quantized algebra of algebraic functions on a group $G$ as a
{subspace}
of $\cU_q(\fg)^*$ {with a Hopf algebra structure}.
In these notes we are taking $\fg=\fsl_2(\C)$, but if we had defined $\cU_q(\fg)$ for
other semisimple Lie algebras $\fg$, the same definition would apply.

\begin{definition}
Let $V$ be a finite dimensional representation of a quantized
enveloping algebra $\cU_q(\fg)$.
We define a \textit{matrix coefficient} of $V$ to be an the element of
the dual $\cU_q(\fg)^*$ of the form
$$
\langle \eta|\cdot |\xi\rangle_V:X \mapsto (\eta, X\xi), \qquad X \in \cU_q(\fg), \, \xi \in V, \, \eta \in V^*
$$
We denote such matrix coefficient as $\langle \eta|\cdot |\xi\rangle_V$.

We define the \textit{quantized algebra of functions} $\cO(G_q)$ as the subspace
of $\cU_q(\fg)^*$ consisting of the matrix elements of all finite-dimensional integrable representations of $\cU_q(\fg)$.
\end{definition}

From the definition, it is not clear why this is an algebra.  In fact, it is not even immediately clear why this is a linear subspace of $\cU_q(\fg)$.  Linearity follows from the fact that 
\[
 \ip{\eta|\cdot|\xi}_V + \ip{\eta'|\cdot|\xi'}_W = \ip{\eta\oplus\eta'|\cdot|\xi\oplus\xi'}_{V\oplus W}.
\]
The other Hopf algebra operations will come from other operations on the $\cU_q(\fg)$-representations $V$ and $W$, as the following theorem shows.

\begin{theorem}
\label{thm:OGq}
Let the notation be as above. Then $\cO(G_q)$ is a Hopf algebra
with multiplication $\mu$ and comultiplication $\Delta$ given as follows:
\begin{equation}
\label{eq:OGq_product}
\langle \eta_1|\cdot |\xi_1\rangle_{V} \langle \eta_2|\cdot |\xi_2\rangle_{W}=
\langle \eta_2\otimes \eta_1|\cdot |\xi_2 \otimes \xi_1\rangle_{V\otimes W},
\end{equation}
\begin{equation}
\label{eq:OGq_coproduct}
\Delta(\langle \eta|\cdot |\xi\rangle_{V})=
\sum_i \langle \eta|\cdot |e_i\rangle_{V}\otimes \langle e^i|\cdot
|\xi\rangle_{V},
\end{equation}
where $\{e_i\}$ and $\{e^i\}$
are bases of $V$ and $V^*$ respectively with $e^i(e_j)=\de_{ij}$.
The unit and counit are given by
\begin{align*}
&1=\langle 1|\cdot | 1\rangle_{\C},  &
&\epsilon(\langle \eta|\cdot |\xi\rangle_{V})=(\eta,\xi),
\end{align*}
where $\C$ denotes the trivial representation of $\cU_q(\lie{g})$ given by the counit of the quantized enveloping algebra $\epsilon:\cU_q(\lie{g})\to \C = \End(\C)$.  We omit the definition of the antipode, which is defined in terms of the contragredient representation, see \cite{KliSch} or \cite[\S{}3.10]{VoiYun:CQG}.
\end{theorem}

We will omit {the proof of} this theorem, although it is not particularly difficult if we take inspiration from the notion of dual Hopf algebras, see Appendix \ref{Hopf-app}.  In fact, the quantum groups $\cU_q(\fg)$ and $\cO(G_q)$ fall into the general philosophy of Pontryagin duality described in Section \ref{sec:Pontryagin}, as explained in the following
{observation, whose details are left to the reader}.

\begin{observation} 
\label{ex:pairing}
Let us denote the pairing of $\cU_q(\fg)$ and $\cO(G_q)$ by
\[
  \cU_q(\fg) \times \cO(G_q) \to \C ; \quad (X,\ip{\eta|\cdot|\xi}_V) = \ip{\eta,X\xi}.
\]
We extend this to a pairing between $\cU_q(\fg)\otimes\cU_q(\fg)$ and $\cO(G_q) \otimes \cO(G_q)$ by $(X\otimes Y, a\otimes b) = (X,a)(Y,b)$.  Then the Hopf algebra operations on $\cU_q(\fg)$ and $\cO(G_q)$ are skew-dual to one another in the sense that, for all $X,Y\in\fg$, $a,b\in\cO(G_q)$,
\begin{align}
\label{eq:skew-pairing}
  (\Delta X, b\otimes a) &= (X,ab) , & (X\otimes Y ,\Delta a) &= (XY,a), \\
  (1_{\cU_q(\fg)},a) &= \counit(a), & (X,1_{\cO(G_q)}) &= \counit(X), \nonumber \\
  (S(X),a) &= (X,S(a)). \nonumber
\end{align}
We have included the relation between the antipodes in this
{list}, although the explicit definition of the antipode on $\cO(G_q)$
was not specified in Theorem \ref{thm:OGq}.   In fact, this last relation can be taken as a definition of $S(a)$ for $a\in\cO(G_q)$ and one can check that it satisfies the antipode relation.
\end{observation}

The pairing
$$
\cU_q(\fg) \times \cO(G_q) \lra \C
$$
is in analogy with the differential action of vector fields on functions
described in \eqref{diff-ac}.  Specifically, if $\fg$ is the Lie algebra of a classical Lie group $G$, and $V$ is a finite-dimensional representation of $V$, then we can alternatively interpret a matrix coefficient $\ip{\eta|\cdot|\xi}_V$ as a function
\[
  \ip{\eta|\cdot|\xi}_V:G \mapsto (\eta,\pi_V(g)\xi),
\]
and then we have $(X,a) = \left.\frac{d}{dt}\right|_{0} a(e^{tx})$ for $a=\ip{\eta|\cdot|\xi}_V\in\cO(G)$ {and $X\in\cU(\fg)$}.  

This defines the algebra of polynomial functions on the quantum group $G_q=\SL_q(2,\C)$, viewed as a complex algebraic variety.

\subsection{Compact Quantum Groups}
In the classical case, the algebra $\cO(G)$ is isomorphic to the $\cO(K)$, where $K$
is the maximal compact subgroup.  This is because any complex polynomial on $G$ is determined by its restriction to $K$, and any polynomial on $K$ extends holomorphically to $G$.  But {as we have seen previously in Section \ref{sec:real_structures},} the polynomials on $K$ admit an additional $*$-operation, given by complex conjugation.  This too has an analogue in the quantum case.

\begin{theorem}
The Hopf algebra $\cO(G_q)$ becomes a Hopf $*$-algebra when equipped with the involution determined by duality with $\UqRk$ as follows:
\[
 (X, a^*) := \overline{(S^{-1}(X)^*,a)},
\]
where the $*$-structure on $\UqRk$ is that given in Proposition \ref{prop:UqRk}.
\end{theorem}

The proof is an exercise using the axioms of a Hopf $*$-algebra.

\medskip

Let us look more carefully at the example of $\fg=\fsl_2(\C)$.

\begin{example}\label{ex-qf}
We specialize to $K=\mathrm{SU}(2)$ and its complexification $G=\mathrm{SL}(2,\C)$.
{The polynomial algebra}, $\cO(G_q)$ is spanned by the matrix coefficients
$\langle \eta|\cdot |\xi\rangle_{V(m)}$, where the representations $V(m)$
are those defined {in Section \ref{sec:Vm} with} $m\in \frac12\N$.   The subspace of matrix coefficients of a given representation $V(m)$ is equal to $(\End (V(m)))^*$ as a matrix coalgebra, meaning it comes with the matrix coproduct \eqref{eq:OGq_coproduct}.  We get
\begin{equation}
\label{eq:Peter-Weyl}
\cO(G_q) = \bigoplus_{m \in \frac12\N} (\End(V(m))^*
\end{equation}
as a direct sum of coalgebras.  The product on $\cO(G_q)$ is more complicated in this picture since it depends upon the decomposition of tensor products into direct sums of irreducibles, see \eqref{eq:OGq_product}.

The representation $V(1)$ is generating for the set of finite dimensional representations of $\mathcal{U}_q(\lie{g})$, in the sense that every $V(m)$ is a subrepresentation of a tensor power $V(1)^{\otimes m}$.  As a consequence, the matrix coefficients for $V(1)$ are algebra generators of $\cO(G_q)$.  

{We remark that} Woronowicz \cite{Woronowicz:pseudogroups} discovered the quantized function algebra $\cO(SL(2,\C))$ by experimentation.  He denoted the generators by
\begin{align*}
  \alpha &= \ip{v^\half|\cdot|v_\half} , &\beta = -q\gamma^* &= \ip{v^\half|\cdot|v_{-\half}}, \\
  \gamma &= \ip{v^{-\half}|\cdot |v_\half}, &\delta = \alpha^* &= \ip{v^{-\half}|\cdot|v_{-\half}},
\end{align*}
where $(v^\half,v^{-\half})$ is the dual basis of $V(1)^*$ to $(v_\half,v_{-\half})\in V(1)$.  These generators satisfy the relations
\begin{align*}
  \alpha\beta&=q\beta\alpha, &  \alpha\gamma &= q\gamma\alpha,  
  &\beta\delta&=q\delta\beta, &\gamma\delta&=q\delta\gamma, &
  \beta\gamma&=\gamma\delta, 
\end{align*}
\begin{align*}  
  \alpha\delta-q\beta\gamma&= \delta\alpha-q\beta\gamma= 1
\end{align*}

For more details see \cite[Ch.4]{KliSch}.

\end{example}

\subsection{{The convolution algebra polynomials on $K_q$}}
As suggested by the discussion of Subsection \ref{sec:Pontryagin}, the quantized enveloping algebra $\UqRk$ is not the only possible model for a convolution algebra dual to the function algebra $\cO(K_q)$.
Indeed, given the decomposition of $\cO(K_q)$ into finite dimensional subspaces, which we saw in Equation \eqref{eq:Peter-Weyl},
it is natural to make the following definition.

\begin{definition}
\label{def:DKq}
We define the algebra of polynomial densities on $\cO(K_q)$ 
as:
\begin{equation}
\label{eq:DKq}
\cD(K_q)=\bigoplus_{m\in \half\N}\End(V(m)).
\end{equation}
The product  is defined entrywise.
\end{definition}

This space $\cD(K_q)$ has a canonical pairing with $\cO(K_q)$ given by applying the natural pairing in each of the components.  Concretely, this amounts to
\[
\End(V(m)) \times \End(V(m))^* \to \C; \qquad (x,\ip{\eta|\,\cdot\,|\xi}_{V(m)} ) \mapsto (\eta, x\xi).
\]
Moreover, $\cD(K_q)$ can be endowed with Hopf $*$-algebra operations by using exactly the same skew-duality relations as we had for $\cU_q(\fg)$ in Example \ref{ex:pairing}, with one technical caveat
{which we mention shortly}.  

For $x,y\in\cD(K_q)$ and $a,b\in\cO(K_q)$, we define
\begin{align}
\label{eq:skew-pairing2}
  (\Delta x, b\otimes a) &= (x,ab) , & (x\otimes y ,\Delta a) &= (xy,a), \\
  (1_{\cD_q(\fg)},a) &= \counit(a), & (x,1_{\cO(G_q)}) &= \counit(x), \nonumber \\
  (S(x),a) &= (x,S(a)), & (x^*,a) &= \overline{(x,S(a)^*)}. \nonumber
\end{align}
The technical caveat is that the coproduct $\Delta$ will have image not in the tensor product $\cD(K_q) \otimes \cD(K_q)$, but in the multiplier algebra of this.
With these operations, $\cD(G_q)$ is a \emph{multiplier Hopf $*$-algebra}.  We will not define multiplier Hopf algebras here, although they are a relatively simple generalisation of Hopf algebras.  For details see \cite{vanDaele:multiplier_Hopf_algebras,vanDaele:duality} or \cite[Ch.2]{VoiYun:CQG}.  

To give a quick idea of why multipliers are necessary, note that the unit for $\cD(K_q)$ is given by
\[
 1= \prod_{m\in \half\N} I_{V(m)},
\]
which belongs in the direct product of the $\End(V(m))$, not the direct sum as in \eqref{eq:DKq}.  Thus $1$ is not an element of the algebra $\cD(K_q)$, but it is a multiplier of $\cD(K_q)$.  

The family of representations of $\UqRk$ on each $V(m)$ yields an embedding 
\[
  \UqRk \hookrightarrow \prod_{m\in\half\N} \End(V(m)).
\]
In this way, the quantized enveloping algebra $\UqRk$ is a $*$-subalgebra of the multiplier algebra of $\cD(K_q)$.  
The algebras $\UqRk$ and $\cD(K_q)$ play very similar roles in the theory of quantized compact semisimple Lie groups, with the main difference being that in a typical $*$-representation on a Hilbert space, elements of $\UqRk$ act as unbounded operators, whereas elements of $\cD(K_q)$ act boundedly.

\section{The Peter Weyl theorem for a compact semisimple quantum group}
\label{sec:Peter-Weyl}

A fundamental problem in unitary representation theory of a group is the following: How does the regular representation of $G$ decompose into irreducible components?  
In the case of a compact group $K$, the solution is given by the famous Peter-Weyl decomposition of $L^2(K)$ in terms   of matrix coefficients of irreducible unitary representations.

This problem also  makes sense for quantum groups, at least if they have a real structure.  In this section, we will describe the analogue of the Peter-Weyl theorem in the case of a quantized compact semisimple Lie group, specifically $\mathrm{SU}_q(2)$.

\subsection{Haar Measure} 

In order to get started, we need the notion of Haar measure.

\begin{definition}
Let $\cA(G)$ be a Hopf algebra.
A \textit{left-invariant integral} on $\cA(G)$ is a linear functional
$\phi:\cA(G)\to \C$ such that for all $a\in\cA(G)$,
\[
  a_{(1)} \phi(a_{(2)}) = \phi(a)1
\]
Similarly, a right-invariant integral satisfies $\phi(a_{(1)})a_{(2)} = \phi(a)1$. 

If moreover $\cA(G)$ is a $*$-Hopf algebra and if $\phi$ is positive definite in the sense that $\ip{a,b} := \phi(a^*b)$ is an inner product on $\cA(G)$, then we say that $\phi$ is a
\textit{left (or right) invariant Haar integral}.
\end{definition}

For instance, if $G$ is a compact semisimple Lie group and $\cA(G)$ is the algebra of polynomial functions, then integration against Haar measure is a bi-invariant Haar integral.

There is a $q$-analogue of this Haar integral on any $q$-deformed compact semisimple Lie group.
As usual we will concentrate only on $K_q = \mathrm{SU}_q(2)$.

\begin{theorem}
\label{thm:Kq-Haar}
The linear functional 
$
\phi: \cO(K_q) \lra \C
$
defined by
\begin{enumerate}
\item $\phi(\langle 1\mid\cdot\mid 1 \rangle_{V(0)})=1$,

\item $\phi(\langle \eta\mid \cdot\mid \xi \rangle_{V(m)})=0$, for all $m\neq 0$ and all $\xi\in V(m)$, $\eta\in V(m)^*$.
\end{enumerate}
is a left- and right-invariant Haar integral on $\cO(K_q)$.
\end{theorem}

The analogous formula works for any compact semisimple quantum group: see \cite[Ch.11]{KliSch}.

Using the left-invariant Haar integral $\phi$, we can define a Hilbert space $\cL^2(K_q)$
as the completion of $\cO(K_q)$ with respect to the inner product 
\[
  \ip{f,g} = \phi(f^*g) \qquad (f,g\in\cO(K_q)).
\]

As in the classical case the Hilbert space $\cL^2(K_q)$ admits two representations of $K_q$, called the left and right regular representations.  To specify these representations, recall from Section \ref{sec:real_structures} that we must take the point of view that a unitary representation of the quantum group $K_q$ is given by a $*$-representation of the quantized enveloping algebra $\UqRk$.

\begin{proposition}
\label{prop:hit}
  There are left and right actions of $\UqRk$ on $\cO(K_q)$ defined by 
    \begin{align*}
 X\hit f &= (X,f_{(2)}) f_{(1)}, & f \hitby X &= (X,f_{(1)})f_{(2)},
\end{align*}
respectively, for $X\in\UqRk$ and $f\in\cO(K_q)$, {where the pairing
(,) is as in Obs. \ref{ex:pairing})}.
\end{proposition}

{Let us define} the left and right regular representations of $\UqRk$ on $\cO(K_q)$ by
  \begin{equation}
  \label{eq:regular_reps}
    \rho(X)f = X \hit f, \qquad \lambda(X)f = f \hitby S(X).
\end{equation}
{Note that in the definition of $\lambda$ we are using the standard trick of turning a right representation into a left representation by using the inverse, or rather the antipode.}
These are $*$-representations in the sense that
  \[
    \ip{\rho(X^*)f,g} = \ip{f,\rho(X)g}, \qquad \ip{\lambda(X^*)f,g} = \ip{f,\lambda(X)g},
  \]
  for all $X\in\UqRk$ and $f,g\in\cO(K_q)$.
{The proof of this, and of Proposition \ref{prop:hit},} is a nice exercise in manipulating the Hopf $*$-algebra axioms and the properties of the Haar integral.  Note that we can also define these representations for the convolution algebra $\cD(K_q)$ by simply replacing $X\in\UqRk$ with $x\in\cD(K_q)$ throughout the above formulas.

\medskip
We can also use the Haar integral to define the Fourier transform.

\begin{theorem}
The linear map
\[
 \cF : \cA(K_q) \to \cD(K_q); \quad a \mapsto \hat{a} := \phi(\,\cdot\, a),
\]
is a linear isomorphism.
\end{theorem}

We call this linear isomorphism the \emph{Fourier transform}.

\medskip
The many properties of the Fourier transform can be found in \cite[Ch.2]{VoiYun:CQG}.  For instance, as in the classical case, Fourier transform intertwines the product in $\cD(K_q)$ with the convolution product in $\cA(K_q)$, which is defined by 
\[
  a*b := \phi(S^{-1}(b_{(1)})a)\, b_{(2)}.
\]

\subsection{The Peter-Weyl Theorem}

{We can now describe the} decomposition of the regular representation  given by the Peter-Weyl Theorem,
or equivalently the Schur Orthogonality Relations, which describes the decomposition of the left and
right regular representation on $L^2(K_q)$
into irreducible components.  

Before we go to its statement, we recall that $L^2(K_q)$ comes equipped
with {a $\UqRk\otimes\UqRk$-representation, given by}
the left and right regular representations:
 \[
    ( X\otimes Y) \cdot f = \lambda(X)\rho(Y)f,  \qquad (X,Y \in\UqRk, ~f\in L^2(K_q)),
\]

Next, let us write $\cL_q^2(V(m))$ for the finite dimensional matrix algebra $\End(V(m))$
equipped with the twisted Hilbert-Schmidt inner product
 \[
  \ip{T_1,T_2}_{\cL^2_q(V(m))} = { \frac{1}{\dim_q V(m)} } \Tr( T_1^* T_2 \pi_m(q^H)) ),
 \]
where  $\dim_q V(m) := \Tr(\pi_m(q^H))$ is called the \emph{$q$-dimension} of the representation $V(m)$.  
 We equip this space with the representation of $\UqRk\otimes \UqRk$ given by 
  \[
    (X\otimes Y)\cdot T = \pi_m(X) T \pi_m(S(Y)).
  \]

\begin{theorem} (Peter-Weyl/Schur Orthogonality)
\label{thm:Peter-Weyl}
Let $K_q=\SU_q(2)$, $q\in(0,\infty)$.  We have an isometric isomorphism of unitary
$K_q\times K_q$-representations:  
 \begin{equation}
 \label{eq:Peter-Weyl}
  L^2(K_q) \cong   \bigoplus_{\scriptstyle m \in \half\N} \cL_q^2(V(m)).
 \end{equation}
  The isometric isomorphism is given by $f\mapsto \bigoplus\pi_m(\hat{f})$, which is to say that for all $f\in\cO(K_q)$ we have
  \[
  \ip{f,g}_{L^2(K_q)} = \sum_{m\in\half\N}   { \frac{1}{\dim_q V(m)} } \Tr(\pi_{m}(\hat{f})^*\pi_{m}(\hat{g}) \pi_m(q^H)),
  \]
  and  this map intertwines the $\UqRk\otimes\UqRk$-representations.
\end{theorem}

The family of functionals $m \mapsto \frac{1}{\dim_q V(m)}\Tr(\,\cdot\, q^H) $ can be understood as the \textit{noncommutative Plancherel measure} on the unitary dual of $K_q$.  The Peter-Weyl Theorem above can be generalized to any compact quantum group, see \cite[Ch.11]{KliSch}.

\section{Quantized complex semisimple Lie groups and their representations}
\label{sec:Gq}

In this section, we will discuss the unitary representation theory of the quantization $G_q$ of a complex semisimple Lie group $G$.  This necessarily means that we will be considering the complex group $G$ as a real Lie group, since unitary representations are not holomorphic maps.

\subsection{Complex semisimple groups as real Lie groups}

We must begin by describing the real structure of the complex semisimple group $G$ in a way that is appropriate for the $q$-deformation process.  This story, in the context appropriate for operator algebras and
unitary representation theory, begins the seminal work of Podle\'{s} and
Woronowicz \cite{PodWor:Lorentz} on the group $\SL_q(2,\C)$, although related ideas had been investigated by Drinfeld
previously \cite{Drinfeld:quantum_groups,Drinfeld:quantum_groups2}.   The key construction is the Drinfeld double, or dually
the Woronowicz double, which were introduced in the above cited articles and which we will describe shortly.

Let us mention that the situtation for general real
semisimple Lie groups is far more delicate.  With the exception of the
quantum group $\SU_q(1,1)$, {the analytic $q$-deformations of noncompact non-complex semisimple Lie groups 
remain very poorly understood}.
For
information about $\SU_q(1,1)$ see \cite{KoeKus:SUq11,DeCommer:SUq11} and references within.  Recent advances can also be found in \cite{DeCDzo:SL2R,DeCDzo:coideals}.

\medskip

Let us consider, therefore, a complex semi-simple Lie group $G$.
To maintain our concrete approach, we will take $G=\SL(2,\C)$.  

Since we are interested in the \emph{unitary} representation theory of $G$ rather than the finite dimensional holomorphic representations which we briefly mentioned in Section \ref{sec:real_structures}, we will not be interested in a quantum analogue of the complex enveloping algebra $\cU(\lie{g})$, but rather in a quantum analogue of the enveloping algebra $\cU^\R(\lie{g}_\C)$ of the \emph{complexified Lie algebra} $\lie{g}_\C = \lie{g}\otimes_\R\C$.  
Since $\lie{g}$ is already a complex Lie algebra, taking its complexification may lead to some confusion.  

As a vector space we have
\[
  \lie{g}_\C \cong \lie{g} \oplus \lie{g}.
\]
It turns out that this is true also as an isomorphism of complex Lie algebras, although the map which institutes the Lie algebra isomorphism is perhaps not immediately obvious.
We will return to this point shortly.
Instead, we will  begin with the philosophy of the quantum double.

The Iwasawa decomposition says that $G$
can be decomposed as a product of real Lie groups
\begin{equation}
\label{eq:KAN}
  G \cong K \bowtie AN,
\end{equation}
where 
\begin{itemize}
\item $K$ is the maximal compact subgroup,
\item $A$ is the subgroup of diagonal matrices with strictly positive entries, 
\item $N$ is the subgroup of unipotent upper triangular matrices $\begin{pmatrix} 1 & z \\ 0 & 1 \end{pmatrix}$.
\end{itemize}
The bowtie in \eqref{eq:KAN} signifies the fact that $K$ and $AN$ are both
Lie subgroups
but neither is a normal subgroup.  Nonetheless, the product map 
\[
K\times AN \to G; \quad (k,an) \mapsto kan
\]
is a diffeomorphism, as is the flipped map $AN \times K \to G$.
As a consequence, we have exchange relations 
\[
  kan = a'n'k',
\]
where $a'n'$ and $k'$ can be calculated as functions of $k$ and $an$.

We want to quantize this setting. 
The philosophy for this is referred to as the quantum duality
principle (see \cite{Gavarini:quantum_duality} and refs therein).
It is not necessary to understand this principle for
what follows, since the quantum group $\SL_q(2,\C)$ will be defined
concretely in terms of generators and relations.  But the principle is
interesting, so we will take a short informal detour to sketch out the idea.

\subsection{The Quantum Duality Principle and the Drinfeld double}

The subgroups $K$ and $AN$ of $G$ are subgroups of equal dimension and both admit natural Poisson Lie-group structures.  A Poisson structure on a Lie group $\lie{k}$ is given by a bivector field $\Pi$.  Linearising $\Pi$ at the identity element gives a linear map $d_e \Pi : \lie{k} \to \lie{k}\wedge \lie{k}$.  Moreover, the axioms of the Poisson bivector turn out to correspond exactly to the fact that $d_e\Pi$ is the dual of a Lie bracket $[\,\cdot\,,\,\cdot\,]$ on the dual space $\lie{k}^*$.  This yields Drinfeld's equivalence of categories between simply connected Poisson Lie groups and Lie bialgebras.  

In the example of $\lie{k}=\lie{su}(2)$, it turns out that the dual Lie algebra structure on $\lie{k}^*$ is isomorphic to $\lie{an}$.  Accordingly, we say that the groups $K$ and $AN$ are \emph{Poisson dual}.

The quantum duality principle states that, upon quantization,
Poisson duality is replaced by Pontryagin duality.  That is, the appropriate quantum analogue of     {the Poisson dual} $AN$ is the Pontryagin dual $\widehat{K}_q$.  In terms of Hopf $*$-algebras, this means
\[
  \cO(AN_q)  = \cO(\widehat{K}_q) = \cD(K_q) \quad \text{and} \quad \cD(AN_q)  =\cD(\widehat{K}_q) = \cO(K_q).
\]

The above discussion justifies the following definition.  As we said above, this definition could also be understood independently of the quantum duality principle.

\begin{definition}
\label{def:DGq}
Let $G=\SL(2,\C)$ (or any other simply connected complex semisimple Lie group).  
We define the space of regular distributions on $G_q$ to be
\[
\cD(G_q) = \cD(K_q \bowtie AN_q) := \cD(K_q) \otimes \cO(K_q).
\]
We equip $\cD(G_q)$ with the twisted product
\begin{equation}
\label{eq:DGq_product}
  (x \bowtie a)(y \bowtie b) = (y_{(1)},a_{(1)})\: x y_{(2)} \bowtie a_{(2)} b \: (S(y_{(3)}),a_{(3)}),
\end{equation}
where we use the Sweedler notation $(1\otimes \Delta)\Delta(a) = a_{(1)} \otimes a_{(2)} \otimes a_{(3)}$, \emph{etc} and we denote with $x \bowtie a$ elements in $\cD(G_q)$.
We also equip $\cD(G_q)$ with the untwisted coproduct 
\begin{equation}
\label{eq:DGq_coproduct}
 \Delta(x\bowtie a) = (x_{(1)} \bowtie a_{(1)}) \otimes   (x_{(2)} \bowtie a_{(2)}).
\end{equation}
\end{definition}

Some remarks on this definition may be helpful.  Firstly, the pairings on the left and right of the right-hand side of Equation \eqref{eq:DGq_product} are just complex numbers, placed on the left and the right for purely aesthetic reasons.  If we use the embeddings
\begin{align*}
  \cD(K_q) &\hookrightarrow \cD(G_q),  & x &\mapsto x\bowtie 1, \\
  \cD(AN_q) = \cO(K_q) &\hookrightarrow \cD(G_q),  & a &\mapsto 1\bowtie a,  
\end{align*}
to identify $\cD(K_q)$ and $\cO(K_q)$ as subalgebras of $\cD(G_q)$, then the product law \eqref{eq:DGq_product} can be seen as an exchange relation:
\begin{equation}
\label{eq:exchange_relation}
    ay= (y_{(1)},a_{(1)})\:  y_{(2)}  a_{(2)}  \: (S(y_{(3)}),a_{(3)}),
\end{equation}
which is analogous to the exchange relation for elements of $K$ and $AN$ in $G$.  Further, these identifications allow us to define an antipode and star operations on $\cD(G_q)$, via
\begin{equation}
\label{eq:Gq_antipode}
   S(xa) = S(a)S(x), \qquad (xa)^* = a^*x^*,
\end{equation}
followed by the exchange relation \eqref{eq:exchange_relation}.  The antipodes and involutions on the right-hand sides of the equations \eqref{eq:Gq_antipode} are those of $\cO(K_q)$ and $\cD(K_q)$, respectively. 

On the other hand, the untwisted coproduct \eqref{eq:DGq_coproduct} corresponds to the classical diffeomorphism $G\cong K\times AN$ as a Cartesian product, since morally, it corresponds to the decomposition $C_0(G) = C(K) \otimes C_0(AN)$ as a tensor product of $C^*$-algebras.  To summarize all of this in the language of noncommutative topology, the quantum group $G_q$ is a Cartesian product $G_q = K_q \times \widehat{K}_q$, but with a twisted group law.

Finally, we define a unit and counit on $\cD(G_q)$ by
\begin{align*}
  &1 = 1\bowtie 1 \in \cM(\cD(K_q)) \otimes \cO(K_q)), \\
  & \epsilon = \epsilon\otimes\epsilon : x\bowtie a \mapsto \epsilon(x)\epsilon(a).
\end{align*}

\begin{theorem}
  The space $\cD(G_q)$ defined above, equipped with the given product, coproduct, unit, counit, antipode and involution  is a multiplier $*$-Hopf algebra.
\end{theorem}

The construction in Definition \ref{def:DGq} can be generalized to any dual pair of Hopf algebras (or multiplier Hopf algebras), and is called the \textit{Drinfeld double}. This is not the construction made by Podle\'{s} and Worowicz in \cite{PodWor:Lorentz}, but rather is the  dual of it.  The Woronowicz double is used to define the algebra of functions
$$
\cO(G_q) := \cO(K_q) \otimes \cD(K_q),
$$
equipped with a twisted coproduct and untwisted pointwise multiplication.  See \cite{PodWor:Lorentz} or \cite[\S{}4.2.4 \& \S{}4.4]{VoiYun:CQG} for more details.

We can also define the quantized enveloping algebra $\cU_q^\R(\lie{sl}_2(\C))$ by
\[
  \UqRg := \UqRk \bowtie \cO(K_q),
\]
with operations given by exactly the same rules as for $\cD(G_q)$ in Definition \ref{def:DGq}, but replacing $x,y\in\cD(K_q)$ with $X,Y\in\UqRk$.  

Note the important appearance of $\R$ in the notation here.   This is to signify that $\UqRg$ is to be thought of as the quantized enveloping algebra of the \emph{complexification} of $\lie{g}=\lie{sl}_2(\C)$, equipped with a real structure, \emph{i.e.} with an involution $*$, which permits us to study the unitary representation theory of $\mathrm{SL}_q(2,\C)$ as described in Section \ref{sec:real_structures}.   The algebra $\UqRg$ should not be confused with the smaller algebra $\cU_q(\lie{g})$, which is just the enveloping algebra $\UqRk$ without its $*$-structure, and is only useful for studying holomorphic representations of $\mathrm{SL}_q(2,\C)$.

Thus $\UqRg$ is a Hopf $*$-algebra, whose elements are multipliers of $\cD(G_q)$.  It follows that modules\footnote{Strictly speaking, we should qualify this by taking \emph{essential} modules for the multiplier Hopf algebra $\cD(G_q)$.  For details, see \cite{VoiYun:CQG}.} for $\cD(G_q)$ are also modules for $\UqRg$, and vice-versa.

\subsection{Representation theory of $\SL_q(2,\C)$}

In the philosophy of quantum groups, unitary representations of $G_q = \SL_q(2,\C)$ correspond to $*$-representations of the convolution algebra 
$\cD(G_q)$.   Since $\cD(G_q) = \cD(K_q) \bowtie \cO(K_q)$, it follows that a representation $\pi$ of $\cD(G_q)$ on a vector space $\cH$ amounts to
\begin{itemize}
\item a representation $\pi$ of $\cD(K_q)$ on $\cH$;
\item a representation $\pi$ of $\cO(K_q)$ on $\cH$
\item a compatibility condition between these representations determined by the exchange relation \eqref{eq:exchange_relation}, namely  for any $x\in\cD(K_q)$, $a\in\cO(K_q)$, 
\begin{equation}
\label{eq:YD_relation}
    \pi(a)\pi(x) = (x_{(1)},a_{(1)})\:  \pi(x_{(2)}) \pi( a_{(2)} ) \: (S(y_{(3)}),a_{(3)}),
\end{equation}
\end{itemize}

\begin{remark}
A representation of $\cD(K_q)$ can equivalently be described via the dual notion of a \emph{corepresentation} of $\cO(K_q)$. 
We will not  {seriously} develop the theory of corepresentations here, for details see \cite[Ch.11]{KliSch},  {but let us make some quick remarks}.  
 {Given a} representation of $\cD(G_q)$ {in the above sense}, one can define a map on $\cH$,
\[
  \check{\pi} : \cH \to  \cO(K_q) \otimes \cH; \qquad \xi \mapsto \xi_{(-1)} \otimes \xi_{(0)} \quad\text{(Sweedler notation)}
\]
which is dual to the representation $\pi$ of $\cD(K_q)$ in the sense that 
\[
\pi(x) \xi = (S(x),\xi_{(-1)}) \xi_{(0)},
\]
for all $x\in \cD(K_q)$, $\xi\in\cH$.
Then the compatibility condition between the representations $\pi$ of $\cD(K_q)$ and $\cO(K_q)$ transforms to the \emph{Yetter-Drinfeld condition}
\[
  (\pi(a)\xi)_{(-1)} \otimes (\pi(a)\xi)_{(0)} = a_{(1)}\xi_{(-1)}S(a_{(3)}) \otimes \pi(a_{(2)})\xi_{(0)}.
\]
For details, see \cite[\S13.1.3]{KliSch}  and \cite[Ch.6.2]{VoiYun:CQG}.
\end{remark}

The following is a crucial example of a $\cD(G_q)$-representation.

Fix a pair of parameters $(\mu,\la) \in \half\Z \times \C$.  
Recall that elements $\mu\in\half\Z$ correspond to integral weights of $K_q$.  Concretely, this means that associated to $\mu$ there is a character $\chi_\mu$ of the Cartan subalgebra 
\[
\cU_q(\lie{h}) = \mathrm{span}\{q^{nH} \mid n\in\Z\} \subset \cU_q(\lie{sl}_2)
\]
given by
$  \chi_\mu:H \mapsto 2\mu$, or more rigorously by
\begin{equation}
\label{eq:chi_mu}
  \chi_\mu : q^H \mapsto q^{2\mu}.
\end{equation}
This character is a $*$-homomorphism, so in fact should be seen as a \emph{unitary} character of $\cU_q^\R(\lie{t}) \subset \cU_q^\R(\lie{su}(2))$, where $\lie{t}$ denotes the Cartan subalgebra of the real Lie algebra $\lie{su}(2)$.  The fact that $\mu\in\half\Z$ means that, classically, $\mu$ integrates to a group character of the diagonal torus subgroup $T \subset \SU(2)$.

With this character, we define the set of $\chi_\mu$-equivariant functions on $K_q$,
\begin{align*}
 \Gamma(\cE_{\mu}) 
  &= \{ f\in \cO(K_q) \mid f_{(1)}  (H,f_{(2)}) = 2\mu f \}.
\end{align*}

\begin{remark}
This is the $q$-analogue of the space of sections of the induced line bundle $\cE_\mu = K \times_T \C_\mu$, on the symmetric space $K/T$,
where $\C_\mu$ denotes the one dimensional representation given by the
character $\chi_\mu$.  
\end{remark}

Put $\cH_{\mu,\la} = \Gamma(\cE_{\mu})$.  We equip space this first with the left-regular representation of $\cD(K_q)$  {as in Equation \eqref{eq:regular_reps}}:
$$
\pi_{\mu,\la}(x)f=f \hitby S(x) = (S(x),f_{(1)}) f_{(2)}.
$$
Next, we add a representation of $\cO(K_q)$ on $\cH_{\mu,\la}$ by the twisted adjoint action:
$$
\pi_{\mu,\la}(a)f=(q^{(\la+1)H},a_{(0)})\, a_{(1)}f S(a_{(3)})
$$

\begin{remark}
Note from this formula that two representations $\pi_{\mu,\la}$ and $\pi_{\mu,\la'}$ are in fact identical if $\la-\la' \in \frac{2\pi i}{\log{q}}\Z = i\hbar^{-1}\Z$, where we put $\hbar = \log(q)/2\pi$.  Therefore, instead of taking our parameters $(\mu,\la) \in \half\Z\times\C$, we may take them to be in $\half\Z \times \lie{h}_q$ where we put
\[
  \lie{h}_q = \C / i\hbar^{-1}\Z.
\] 
\end{remark}

The following theorem states that the above representations of $\cD(K_q)$ and $\cO(K_q)$ do indeed define a representation of $\SL_q(2,\C)$ on $\cH_{\mu\,\la}$, via
\[
  \pi_{\mu,\la}(x\bowtie a) = \pi_{\mu,\la}(x)\pi_{\mu,\la}(a), \qquad (x\in\cD(K_q),~a\in\cO(K_q)).
\]  
It is called the \emph{principal series representation} of parameter $(\mu,\la) \in \half\Z \times \lie{h}_q$.  

\begin{theorem}
\label{thm:quantum_principal_series}
Let the notation be as above.
\begin{enumerate}
\item The representations $\pi_{\mu,\la}$ of $\cD(K_q)$ and $\cO(K_q)$ are compatible in the sense of Equation \eqref{eq:YD_relation}.  Thus, $\cH_{\mu,\la}$ is a representation of $\cD(G_q)$, although not necessarily a $*$-representation.  
\item The representations extend continuously to the norm closure $H_{\mu,\la}$ of $\cH_{\mu,\la}$ in $L^2(K_q)$.
\item The principal series representation $\pi_{\mu,\la}$ is unitary, \emph{i.e.},
it  is a $*$-representation of $\cD(G_q)$, if and only if the parameter $\la$ is purely imaginary, meaning $\la \in i\lie{t}_q$ where $\lie{t}_q = \R / \hbar^{-1}\Z$.

\item The representation $\pi_{\mu,\la}$ is irreducible if and only if
$$
\pm \la \not\in |\mu|+i\pi h^{-1}\Z+\N^\times
$$
In particular, the unitary representations are always irreducible. 

\item Two irreducible principal series representations $H_{\mu,\la}$ and $H_{\mu',\la'}$ are equivalent if and only if $(\mu,\la) = \pm(\mu,\la)$.
\end{enumerate}
\end{theorem}

 {
\begin{remark}
Item 3 in the above theorem concerns only the representation on $\cH_{\mu,\la}$ with its inner product inherited from $L^2(K_q)$.  It may be the case that the representation on $\cH_{\mu,\la}$ becomes unitary when equipped with a different inner product.  This is the case for the \emph{complementary series} of representations. 
The complementary series representations will not play a role in what follows, so we will ignore them here.  For details, see \cite[\S6.10]{VoiYun:CQG}
\end{remark}
}
As elsewhere, these results all have analogues for a general complex semisimple Lie group $G$.  The parameter space becomes $\mathbf{P} \times \lie{h}_q$, where $\mathbf{P}$ is the integral weight lattice and $\lie{h}_q = \lie{h}/i\hbar^{-1}\mathbf{Q}^\vee$ with $\mathbf{Q}^\vee$ being the coroot lattice.  The condition for equivalence in the final point is replaced by requiring the pairs $(\mu,\la)$ and $(\mu',\la')$ to be in the same orbit of the Weyl group on the parameter space.

For proofs of these facts, see Ch.6 of \cite{VoiYun:CQG}.  The case of $\SL_q(2,\C)$ is dealt with in Sections 6.7 and 6.10.

\section{The Plancherel formula for a complex semisimple quantum group}

 {
In Section \ref{sec:Peter-Weyl}, we discussed the problem of decomposing the regular representation of the compact quantum group $K_q$ into a direct sum of irreducibles.  The same question can be posed for the decomposition of the complex semisimple quantum group $G_q$.  In this case, and in the classical counterpart for $G$, the decomposition of $L^2(G)$ is no longer a direct sum but a direct integral.  This is called the problem of the \emph{Plancherel measure} because the solution relies on finding a measure on the set of irreducible unitary representations of $G$ or $G_q$.  The result generalizes the Plancherel Theorem for the group $\R$, stating that Fourier transform is an isometric isomorphism on $L^2(\R)$  with respect to an appropriate rescaling of the Lebesgue measure.
}

 {
In the classical case, the Plancherel measure for complex semisimple Lie groups was found by Harish-Chandra.  In this section, we will discuss the analogue of Harish-Chandra's theorem for the quantized complex semisimple Lie groups, particularly $\mathrm{SL}_q(2,\C)$.  
}

\subsection{The classical Plancherel formula for $\SL(2,\C)$}

Recall that a Lie group $G$ is called \textit{unimodular}, if the right and left invariant Haar measures coincide.
This is a property
shared by all semisimple complex groups.

We define an inner product $\ip{ \,\cdot\,,\,\cdot\,}$ on $\Cc^\infty(G)$ as usual by 
\[
  \ip{f,g} = \int_{x\in G} \overline{f(x)} g(x) \, dx
\]
so that the completion of $\Cc^\infty(G)$ with respect to this measure is $L^2(G)$.

The group $\SL(2,\C)$ has a family of irreducible unitary representations $H_{\mu,\la}$, called the unitary principal series representations, which are induced from unitary characters the Borel subgroup $B=TAN$ of invertible upper triangular matrices.   {They are} indexed by $(\mu,\la)\in \half\Z \times \lie{a}^*$ in a manner   {analogous to the $G_q$-representations in Theorem \ref{thm:quantum_principal_series}}.   We omit the definition here, since our real interest is the quantum group analogue.  Recall also that any unitary representation $\pi:G\to U(H)$ can be integrated to a representation of the convolution algebra $\Cc^\infty(G)$ by the formula
\begin{equation}
\label{eq:integrated_form}
  \pi(u) = \int_{g \in G} u(x) \pi(g)  \, dg  \in \cB(H).
\end{equation}

We also recall that the Hilbert space $\cL^2(H)$ of Hilbert-Schmidt operators is defined as the set of operators {on $H$}of finite norm with respect to the inner product
\[
  \ip{S,T} = \Tr(S^*T).
\]
It is a fact that for every $f\in\Cc^\infty(G)$, the integrated principal series representation $\pi_{\mu,\la}(f)$ is a Hilbert-Schmidt operator, for every $(\mu,\la)\in\half\Z\times\lie{a}^*$.

\begin{theorem}[Plancherel Theorem for $\SL(2,\C)$]
\label{thm:SL2C-Plancherel}
Let $G=\SL(2,\C)$.  
The space $L^2(G)$
decomposes as a direct integral of $G\times G$ representations
\begin{equation}
\label{eq:SL2C-Plancherel}
  L^2(G) \cong \directint_{(\mu,\la)\in\half\Z \times \lie{a}^*} \cL^2(H_{\mu,\la}) \,dm(\mu,\la),
\end{equation}
where the Plancherel measure is
\[
dm = \half |\mu+i\lambda|^2 \,d\mu\,d\la,
\]
with $d\mu$ being counting measure on $\half\Z$ and $d\la$ Lebesgue measure on $\lie{a}^* \cong \R$.  
More explicitly, if $f,g\in\Cc^\infty(G)$, then 
$$
\langle f, g \rangle_{L^2(G)}=
\int_{(\mu,\la)} \Tr(\pi_{\mu,\la}(f)^* \pi_{\mu,\la}(g)) \, dm(\mu,\la).
$$
The isomorphism \eqref{eq:SL2C-Plancherel} intertwines the $G\times G$-representations, where $(g,h)\in G\times G$ acts
\begin{itemize}
\item on $L^2(G)$ by the left and right regular representations $\lambda(g)\rho(h)$ and 
\item on $\cL^2(H_{\mu,\la})$ by $T\mapsto \pi_{\mu,\la}(g)T\pi_{\mu,\la}(h^{-1})$.
\end{itemize}
\end{theorem}

This is the basic case of a more general Plancherel Theorem for all semisimple Lie groups---one of the crowning achievements of Harish-Chandra.

 {
\subsection{The Plancherel formula for $\mathrm{SL}_q(2,\C)$}
}
In 1999, Buffenoir and Roche \cite{BufRoc} proved an analogue of this formula for the quantum group $\SL_q(2,\C)$, which has more recently been generalized to the $q$-deformations of all complex semisimple Lie groups \cite{VoiYun:Plancherel}.  As previously mentioned, we do not currently know how to generalize this to other non compact real semisimple groups, because we do not know how to construct their $q$-deformations.

We shall see that the $q$-deformed Plancherel formula is remarkably similar to Harish-Chandra's formula in Theorem \ref{thm:SL2C-Plancherel} above.  To state it, we need the Haar integral on $\SL_q(2,\C)$.  Let us write $G=\SL(2,\C)$ and $K=\SU(2)$, as before.

We have already defined Haar measure $\phi$ on the compact quantum group $\SU_q(2)$, see Theorem \ref{thm:Kq-Haar}.  {We now define Haar measure on $AN_q = \widehat{K}_q$ and on $G_q$.}

\begin{theorem}
\begin{enumerate}
\item
The quantum group
$$
\cO(\widehat{K}_q) =\cD(K_q) \cong \bigoplus_{m\in\half\N} \End V(m),
$$
which is dual to $\cO(K_q)$, has a right-invariant Haar integral $\hat\psi$ defined by
\[
  \hat\psi : x=\bigoplus_m x_m \mapsto \sum_{m\in\half\N} \dim_q V(m) \Tr( x_m \pi_m(q^{-H})).
\]
\item
The quantum group $\cO(G_q) = \cO(K_q) \otimes \cO(\widehat{K}_q)$ which is dual to $\cD(G_q)$ has a two-sided Haar integral $\Phi$ defined by
\[
 \Phi(a \otimes x) = \phi(a)\hat\psi(x).
\]
\end{enumerate}
\end{theorem}

As usual, the Haar integral $\Phi$ lets us define an inner product on $\cO(G_q)$ by
\[
  \ip{f,g} = \Phi(f^*g),
\]
and we denote the completion by $L^2(G_q)$.  We also have the Fourier transform,
\begin{align*}
  \cF : \cO(G_q) &\to \cD(G_q),&
    f &\mapsto \hat{f}=\Phi(\,\cdot\,f).
\end{align*}
If $\pi$ is a $*$-representation of $\cD(G_q)$ and $f\in\cO(G_q)$, then the analogue of the classical integrated representation  \eqref{eq:integrated_form} is the operator $\pi(\hat{f})$.  In fact, this would be a more accurate notation also for the classical notion of \eqref{eq:integrated_form}, since $f\mapsto \pi(f)$ defines a representation of the convolution algebra of smooth densities not the pointwise algebra of smooth functions.

Finally, we need a twisted variant of the Hilbert-Schmidt operators, similar to that which appeared in the Peter-Weyl formula for $K_q$, Theorem \ref{thm:Peter-Weyl}.  Suppose $\pi$ is a $*$-representation of $\cD(G_q)$ on a Hilbert space $H$ whose restriction to $\cD(K_q)$ is integrable, so that $\pi(q^H)$ is a positive unbounded operator.  We define $\cL^2_q(H)$ to be the set of operators on $H$ of finite norm with respect to the inner product
\[
 \ip{S,T} = \Tr(ST\pi(q^{-H})).
\]

\begin{remark}
The unbounded operator $\pi(q^{-H})$ in this formula is called the \emph{Duflo-Moore} operator.  Hilbert-Schmidt norms twisted by a Duflo-Moore operator can appear in the Plancherel formula for classical groups if the group is not unimodular.  Here, the quantum group $\SL_q(2,\C)$ is unimodular, but a Duflo-Moore operator is nonetheless necessary. 
\end{remark}

\begin{theorem}[Plancherel Theorem for $\SL_q(2,\C)$]
Let $G_q=\SL_q(2,\C)$.  
The space $L^2(G_q)$ decomposes as a direct integral of $\cD(G_q)\otimes \cD(G_q)$\--repre\-sentations
\begin{equation}
\label{eq:SLq2C-Plancherel}
  L^2(G_q) \cong \directint_{(\mu,\la)\in\half\Z \times i\lie{t}_q} \cL_q^2(H_{\mu,\la}) \,dm_q(\mu,\la),
\end{equation}
where the Plancherel measure is
\[
dm_q = \half |[\mu+i\lambda]_q|^2 \,d\mu\,d\la,
\]
with $d\mu$ being counting measure {on $\half\Z$} and $d\la$ being Lebesgue measure on the circle $\lie{t}_q \cong \R/\hbar^{-1}\Z$.  
More explicitly, if $f,g\in\Cc^\infty(G)$, then 
\begin{equation}
\label{eq:Gq-Plancherel}
\langle f, g \rangle_{L^2(G)}=
\int_{(\mu,\la)} \Tr(\pi_{\mu,\la}(\hat{f})^* \pi_{\mu,\la}(\hat{g}) \pi_{\mu,\la}(q^{-H})) \, dm_q(\mu,\la).
\end{equation}

The isomorphism \eqref{eq:SLq2C-Plancherel} intertwines the $\cD(G_q)\otimes \cD(G_q)$-representations, where $u\otimes v \in \cD(G_q)\otimes\cD(G_q)$ acts
\begin{itemize}
\item on $L^2(G)$ by the left and right regular representations, $(u \otimes v) : f \mapsto (S(u),f_{(1)}) f_{(2)} (v,f_{(3)})$ and 
\item on $\cL_q^2(H_{\mu,\la})$ by $(u \otimes v): T\mapsto \pi_{\mu,\la}(u)T\pi_{\mu,\la}(S(v))$.
\end{itemize}

\end{theorem}

\subsection{Proof of the Plancherel Formula, part I}

We will finish these notes with a very rough outline of the proof of the Plancherel Theorem above.
The first thing to note is that, under the Fourier transform, the Haar integral of $\cO(G_q)$ transforms to the counit of $\cD(G_q)$, since for $f\in\cO(G_q)$, 
\[
  \Phi(f) = \hat{f}(1) = \epsilon(\hat{f}).
\]
Therefore, if we put $u=\hat{f}^*\hat{g} \in \cD(G_q)$, the desired Plancherel formula \eqref{eq:Gq-Plancherel} becomes
\begin{equation}
\label{eq:Plancherel_trace}
  \epsilon (u) = \int_{(\mu,\la)\in\half\Z\times i\lie{t}_q} \Tr(\pi_{\mu,\la}(u) \pi_{\mu,\la}(q^{-H})) \, dm_q(\mu,\la).
\end{equation}

It now suffices to consider the case 
\begin{equation}
\label{eq:u-special}
  u = | v_i\rangle\langle v^j | \otimes \ip{v^k| \,\cdot\, | v_l} \in \End V(m) \otimes (\End V(m'))^*,
\end{equation}
where $v_i$ is a standard basis vector for $V(m)$, $v^j$ a dual basis vector, and similarly $v_l\in V(m')$ and $v^k\in V(m')^*$.  In this notation, $| v_i\rangle\langle v^j|$ denotes the rank-one operator $v\mapsto \ip{v^j|v}v_i$ belonging to $\End(V(m))$.  

With such an element $u$, an explicit calculation can be carried out to simplify the right-hand side of \eqref{eq:Plancherel_trace}.
We will omit this calculation and simply state that the result comes to
\begin{align}
\label{eq:Lefschetz}
 \int_{(\mu,\la)\in\half\Z\times i\lie{t}_q} \Tr(\pi_{\mu,\la}(u) &\pi_{\mu,\la}(q^{-H})) \, dm_q(\mu,\la) \nonumber \\
& =   \Tr(\pi_{0,-1}(\tilde{u})) - \Tr(\pi_{1,0}(\tilde{u})),
\end{align}
where
{
$\pi_{0,-1}$ and $\pi_{1,0}$ are the principal series representations of Theorem \ref{thm:quantum_principal_series} and
} $\tilde{u}$ is some explicit element of $\cD(G_q)$ related to $u$ via the antipodes of $\cO(K_q)$ and $\cD(K_q)$.  For the details of this calculation, see Section 6 of \cite{VoiYun:Plancherel}.

\begin{remark}
The analogous calculation can be carried out for a general complex semisimple quantum group $G_q$.  The result comes to
\begin{equation}
\label{eq:Lefschetz2}
  \sum_{w\in W} (-1)^{|w|} \Tr(\pi_{-w.0, -w.0 -2\rho}(\tilde{u})),
\end{equation}
where $W$ is the Weyl group of $\lie{g}$, $|w|$ denotes the Bruhat length of $w\in W$,  $\rho$ is the half-sum of the positive roots and
$w.\lambda = w(\lambda+\rho) -\rho$ is the $\rho$-shifted Weyl group action.
\end{remark}

The alternating sum of traces in \eqref{eq:Lefschetz} and \eqref{eq:Lefschetz2} is reminiscent of the Lefschetz trace formula, and indeed it is this that will lead to the equality \eqref{eq:Plancherel_trace}, passing via the Bernstein-Gelfand-Gelfand complex.  This is the remaining piece of the puzzle.

\subsection{Harish-Chandra bimodules}
\label{HC-bimodules}

Our goal is to convert the study of unitary representations, which is rather difficult,
into the study of Verma modules, which is more algebraic and in principle simpler.  

Consider first the case of the classical group $G=\SL(2,\C)$.  As already mentioned several times---see Sections \ref{sec:real_structures} and \ref{sec:Gq}---we begin by linearizing and complexifying the problem, thus obtaining $\C$-linear representations of the complexified Lie algebra, $\lie{g}_\C = \lie{g}\otimes_\R\C$, where $\lie{g}$ is complex to begin with.  

At this point, we observe that we have an isomorphism of complex Lie algebras,
\begin{equation}
\label{eq:g-complexification}
  \lie{g}_\C \cong \lie{g} \oplus \lie{g}^\op; \qquad 
  X \mapsto (X , \overline{X}^\rmt) \quad \text{(for $X \in \lie{g}$)},
\end{equation}
where $\overline{X}^\rmt$ is the conjugate transpose of $X$.  This in turn induces an algebra isomorphism
\begin{equation}
\label{eq:Ug-complexification}
  \cU(\lie{g}_\C) \cong \cU(\lie{g}) \otimes \cU(\lie{g})^\op.
\end{equation}
As a consequence, the usual procedure of converting a unitary representation $H$ of $G$ into a complex-linear representation of the complexified enveloping algebra $\cU(\lie{g}_\C)$ can be further modified by using \eqref{eq:Ug-complexification} to obtain a $\cU(\lie{g})$-bimodule structure on $H$.  
More precisely, since Lie algebra elements act as unbounded operators, the associated $\cU(\lie{g})$-bimodule will be given by the dense subspace $\cH$ of $K$-finite vectors\footnote{A vector $v$ in $H$ is \emph{$K$-finite} if $\pi(K)v$ is a finite dimensional subspace.  For such vectors, the map $g\mapsto\pi(g)v$ is always smooth, so the derived representation of $\cU(\lie{g}_\C)$ is well-defined.} in $H$.  

{Next consider the embedding of $\lie{k}$ into $\lie{g}$, corresponding to the inclusion of the maximal compact subgroup $K$ in $G$.  {Note that elements $X\in\lie{k}=\lie{su}(2)$ are skew-adjoint, $\overline{X}^\rmt = -X$.  Therefore,} upon complexification we obtain an embedding of $\lie{k}_\C$ into $\lie{g}_\C\cong \lie{g}\oplus\lie{g}^\op$ which, according to}
Equation \eqref{eq:g-complexification}, is given by
\begin{equation}
\label{eq:k-embedding}
  X \mapsto (X,-X)  \quad \text{(for $X \in \lie{k}$)}.
\end{equation}
Thus, the action of $X\in\lie{k}_\C$ on the $\cU(\lie{g})$-bimodule $\cH$ is given by
 \begin{equation}
 \label{eq:k-adjoint}
   v \mapsto X\!\cdot\! v - v\!\cdot\! X.
\end{equation}

 We can extend the action \eqref{eq:k-adjoint} by universality to $\cU(\lie{k}_\C)$.  This motivates the following notation.
 
 \begin{definition}
 \label{def:adjoint_action}
   Let $\cH$ be a $\cU(\lie{g})$-bimodule.  For $X\in \cU(\lie{k}_\C)$, and $v\in\cH$, we write
   \[
     X\hit v = X_{(1)}\cdot v \cdot S(X_{(2)}),
   \]
   and call this the \emph{adjoint action} of $\cU(\lie{g})$ on $\cH$.
 \end{definition}

 Note that this definition may also cause confusion, since $\cU(\lie{k}_\C) \cong \cU(\lie{g})$
 as an algebra.  Therefore the bimodule $\cH$ is now equipped with three actions of $\cU(\lie{g})$:
 \begin{itemize}
 \item a left action, denoted by $X\!\cdot\! v$,
 \item a right action, denoted by $v\!\cdot\!X$,  
 \item an adjoint action of $\cU(\lie{k}_\C)\cong\cU(\lie{g})$, denoted by $X\hit v$,
 \end{itemize}
 all of which come naturally from the unitary representation of $G$ on $H$, as well as its restriction to the maximal compact $K$.   These three actions are compatible in the sense that
 \begin{equation}
 \label{eq:HC-compatibility}
   X\!\cdot\!v = (X_{(1)}\hit v)\cdot X_{(2)},
 \end{equation}
 as one can readily check from Definition \ref{def:adjoint_action}. 
 
 As a technical point, note also that since $K$ is compact, any unitary representation of $K$ decomposes into a direct sum of finite dimensional representations.  Therefore, the adjoint action of $\cU(\lie{k})$ is a \emph{locally finite} representation, meaning that for any $v\in \cH$,  its image $\cU(\lie{k}_\C)\hit v$ is a finite dimensional subspace of $\cH$.  This is typically not true of the other two actions of $\cU(\lie{g})$.

 To summarize all of the above manipulations, any unitary representation $H$ of $G$ gives rise to a $\cU(\lie{g})$-bimodule $\cH$ equipped in addition with a locally finite ``adjoint'' action of $\cU(\lie{g})$ which is compatible in the sense of Equation \eqref{eq:HC-compatibility}.  Such objects, with a few added technical considerations, are called \emph{Harish-Chandra bimodules}.
 
 \medskip
 
Amazingly, this whole procedure can be $q$-deformed.  This is the result of deep work due to Joseph and Letzer, see \cite{Joseph:book}.
One notable technicality is that the algebra $\cU_q(\lie{g})$ will occasionally have to be replaced by its locally finite part $F\cU_q(\lie{g})$ for the adjoint action, whereas in the classical case $F\cU(\lie{g})$ and $\cU(\lie{g})$ coincide.    For the non-expert, this technicality can be ignored, but we will keep it in our notation for accuracy.
 
 The $q$-analogue of the isomorphism {$\lie{g}_\C \cong \lie{g}\oplus\lie{g}$, or its extension to the enveloping algebra} \eqref{eq:Ug-complexification}, is as follows.
\begin{theorem}
\label{thm:iota}
There is an algebra morphism 
\[
\iota : \UqRg := \UqRk \bowtie \cO(K_q)
   \hookrightarrow \cU_q(\lie{g}) \otimes \cU_q(\lie{g})^\op,
\]
with image
\[
((\id\otimes {S})\Delta \cU_q(\lie{g})) .
   ( F\cU_q(\lie{g}) \otimes 1)
\]
or equivalently,
\[
((\id\otimes {S})\Delta \cU_q(\lie{g})) .
   (1\otimes F\cU_q(\lie{g})).
\]
\end{theorem}

To explain this theorem, the first component $(\id\otimes {S})\Delta U_q(\lie{g})$ in the image corresponds to the embedding of the compact part $\lie{k}_\C$ in $\lie{g}\oplus\lie{g}^\op$ that we saw in Equation \eqref{eq:k-embedding}, and indeed the morphism $\iota$ restricted to the compact part $\UqRk\subseteq \UqRg$ is given by
\[
  \iota:X \mapsto (\id\otimes S)\Delta(X)
\]

The other components, with $F\cU_q(\lie{g})$ in the first or second leg, correspond to the components $\lie{g}$ and $\lie{g}^\op$ that we saw in Equation \eqref{eq:g-complexification}.   The action of $\iota$ on $\cO(K_q)$ is given in terms of the {so-called} $l$-functionals, but we will not present the details here.

Taking inspiration from the classical case above suggests that the following definition will be useful.  The final conditions in this definition are technical points which we will not develop.

\begin{definition}
A \emph{quantum Harish-Chandra bimodule} is an $F\cU_q(\lie{g})$-bimodule $\cH$ equipped with an \emph{adjoint action} $\hit$ of $\cU_q(\lie{g})$,  which is compatible with the bimodule structure in the following sense
\begin{equation}
\label{eq:qHC-sompatibility}
 X\cdot v = (X_{(1)}\to v)\cdot X_{(2)},
\end{equation}
for all $X\in\cU_q(\lie{g})$ and $v\in\cH$, and such that
\begin{enumerate}
 \item the adjoint action is locally finite,
 \item the right $F\cU_q(\lie{g})$-action is finitely generated,
 \item the right $F\cU_q(\lie{g})$-action has annihilator of finite codimension.
\end{enumerate}
\end{definition}

Two examples will be of crucial importance.

\begin{example}
 Let $H$ be any $*$-representation of $\cD(G_q)$ and hence of $\UqRg$.  Let $\cH\subseteq H$ be the locally finite part for the action of the subalgebra $\UqRk$.
 Then $\cH$ inherits the structure of a quantum Harish-Chandra bimodule via the morphism $\iota$.

  In particular the principal series representations give rise to quantum Harish-Chandra bimodule structures on $\cH_{\mu,\la}$.
   
\end{example}

\begin{example}
   Let $M,N$ be Verma modules for $U_q(\lie{g})$ as in Observation \ref{obs:Verma}, or more generally let $M,N$ be modules in category $\cO$ (see \cite{Joseph:book, VoiYun:CQG} for the definition of quantum category $\cO$).  We equip the space $\Hom(M,N)$ with 
  \begin{itemize}
   \item left and right actions of $Y,Z \in F\cU_q(\lie{g})$ given by 
   \[
   (Y\cdot\phi\cdot Z)(m) = Y\cdot(\phi(Z\cdot m)),
   \]
   for $\phi\in\Hom(M,N)$, $m\in M$;
   \item an adjoint action of $X\in \UqRk \cong \cU_q(\lie{g})$ given by 
   \[
   (X\hit\phi)(m) = X_{(1)}\cdot(\phi(\hat{S}(X_{(2)})\cdot m).
   \]
  \end{itemize}
  Now put $\cH = F\Hom(M,N)$, the locally finite part of $\Hom(M,N)$ with respect to the adjoint action.  Then $\cH$ is a quantum Harish-Chandra bimodule.  At least for the compatibility condition \eqref{eq:qHC-sompatibility}, this is an easy check.
 \end{example}

Inspired by similar results for classical groups, Joseph and Letzter observed that the above two examples are isomorphic in many cases.  In particular, as the next theorem shows, the unitary principal series representations can be realized as $F\Hom(M,N)$ for appropriate modules $M$ and $N$.  We recall from observation \ref{obs:Verma} that we defined the Verma module
$M(m)$ with any highest weight $m\in\C$.  There is also a notion of dual Verma module $M(m)^\vee$, which we will not detail here, see \cite[\S5.1.2]{VoiYun:CQG}.

  \begin{theorem}[Joseph-Letzter]
  Let $(\mu,\la) \in \half\Z\times i\lie{t}_q$.  Then, as quantum Harish-Chandra bimodules, 
   \[
      \cH_{\mu,\lambda} \cong F\Hom(M(\ell), M(r)^\vee)
   \]
   where $\mu = \ell-r$, $\lambda =-\ell-r-1$.
  \end{theorem}

In fact, with some caveats, Joseph and Letzter's isomorphism defines an equivalence of categories between certain subcategories of category $\cO$ and of quantum Harish-Chandra bimodules.   To be specific, let us fix $\ell=0$.  Let $\cO_0$ denote the subcategory of category $\cO$ in which all weights appearing are integral.  In particular, $\cO_0$ contains the irreducible integrable $\cU_q(\lie{g})$-modules $V(m)$, as well as all of the Verma modules $M(m)$ appearing in the projective resolution \eqref{eq:BGG1} of $V(m)$, and also their dual modules $M(m)^\vee$.  

Joseph and Letzter prove that the map
\begin{equation}
\label{eq:F0}
   N \mapsto F\Hom(M(\ell),N)
\end{equation}
is an exact functor from $\cO_0$ to the category of quantum Harish-Chandra bimodules.  
This means that if we apply the functor \eqref{eq:F0} to {dual of} the resolution  \eqref{eq:BGG1} of the trivial representation
\[
  0 \to  V(0) \to M(0)^\vee \to M(-1)^\vee \to 0,
\]
then we obtain a short-exact sequence of quantum Harish-Chandra bimodules, and hence a resolution of the trivial representation of $\UqRg$ by principal series representations:
\begin{equation}
\label{eq:BGG2}
  0 \to \C \to  \cH_{0,-1} \to \cH_{1,0} \to 0.
\end{equation}

\begin{remark}
\label{rem:BGG}
This resolution of the trivial representation by principal series representations generalizes to all quantized complex semisimple groups.  It is called the \emph{geometric Bernstein-Gelfand-Gelfand resolution} and is a quantum analogue of a well-studied differential complex on flag varieties,
see for instance \cite{CapSloSou:BGG}.  

The BGG complex has seen various applications in the noncommutative geometry of quantum groups, see for instance \cite{HecKol:differential_forms, VoiYun:SUq3,OBuSom:root_vectors}.
\end{remark}

\subsection{Proof of the Plancherel Formula, part II}

Now we can complete the proof of the Plancherel formula.  The principal series representations appearing in the geometric BGG resolution \eqref{eq:BGG2} are precisely those which appear in the reduction of the Plancherel formula, Equation \eqref{eq:Lefschetz}.  Since the maps are morphisms of $\UqRg$-representations, and hence $\cD(G_q)$-representations, they intertwine the action of $\tilde{u}$ on each.  The Lefschetz trace principle implies that the alternating sum of the traces is zero:
\[
  \epsilon(\tilde{u}) - \Tr(\pi_{0,-1}(\tilde{u})) + \Tr(\pi_{1,0}(\tilde{u})) =0.
\] 
From the definition of $\tilde{u}$, which we have elided
(see eq \eqref{eq:Lefschetz}),
one easily obtains $\epsilon(\tilde{u}) = \epsilon(u)$.  Therefore, the integral \eqref{eq:Lefschetz} is equal to $\epsilon(u)$.  This proves the formula \eqref{eq:Plancherel_trace} and hence the Plancherel formula.  

\medskip

The above proof works equally well for the $q$-deformations of all complex semisimple Lie groups.  The projective resolution \eqref{eq:BGG1} of the finite dimensional integral module $V(m)$ by Verma module has a generalization in higher rank, called the \emph{algebraic Bernstein-Gelfand-Gelfand resolution} \cite{BerGelGel}.   The Verma modules which appear in the BGG resolution are those which have highest weights in the $\rho$-shifted Weyl group orbit of $m$:
\[
0 \to \cdots \to \bigoplus_{\substack{w\in W \\ \ell(w)=k}} M(w.m) \to \cdots \to M(m) \to V(m) \to 0.
\]
In the quantum case, this result is due to Heckenberger and Kolb \cite{HecKol:BGG}.

We can then take the dual of this resolution in category $\cO$, and apply Joseph and Letzter's functor as Equation \eqref{eq:BGG2}.  This yields the \emph{geometric Bernstein-Gelfand-Gelfand complex} which we mentioned in Remark \ref{rem:BGG}.    The principal series representations which arise in the geometric Bernstein-Gelfand-Gelfand resolution of the trivial representation are again precisely those that appear in the integral calculation \eqref{eq:Lefschetz2}.  Therefore an application of the Lefschetz trace principal once again completes the proof.


\appendix
\section{Hopf Algebras}\label{Hopf-app}

In this section we recap some key definitions regarding Hopf algebras.
For more details see \cite{FioLle, Montgomery:Hopf_lectures}.

Let $k=\R,\C$ be our ground field. 

\begin{definition} \label{hopfalgebra}
We say that $A$  is a \textit{Hopf algebra} if it
has the following properties:
\begin{enumerate}
\item $A$ is an \textit{algebra} (not necessarily commutative), that is,
  there are linear maps, the {\it  multiplication} $\mu:A \otimes A \lra A$
and the {\it unit} $i:k \lra A$ such that the following diagrams commute
$$
\begin{array}{c}
\begin{CD}
A \otimes k  @> \id \otimes i >> A \otimes A \\
@V{\cong}VV @VV\mu V \\
A  @>\id  >> A
\end{CD}
\qquad
\begin{CD}
k\otimes A  @> i \otimes \id >> A \otimes A \\
@V\cong VV @VV\mu V \\
A  @>\id  >> A
\end{CD}
\\ \\
\begin{CD}
A \otimes A \otimes A @> \mu \otimes \id >> A \otimes A \\
@V{\id \otimes \mu}VV @VV\mu V\\
A \otimes A @> \mu >>  A
\end{CD}
\end{array}
$$
A {\it morphism} $\phi:A \lra B$ of two algebras,  with multiplication $\mu_A$ and $\mu_B$ and unit $i_A$ and $i_B$ respectively, is a linear map such that
$$\mu_B\circ (\phi\otimes \phi)=\phi\circ \mu_A,\qquad \phi\circ i_A=i_B\,.$$

\item  $A$ is a \textit{coalgebra}, that is, we can define
two linear maps called
\textit{comultiplication} $\Delta:A \lra A \otimes A$
and \textit{counit} $\ep: A \lra k$ with the following properties:
$$
\begin{array}{c}
\begin{CD}
A \otimes A  @> \id \otimes \ep >> A \otimes k \\
@A{\Delta}AA @AA\cong A \\
A  @>\id  >> A
\end{CD}
\qquad
\begin{CD}
A \otimes A  @> \ep \otimes \id >> k \otimes A \\
@A{\Delta}AA @AA\cong A \\
A  @>\id  >> A
\end{CD}
\\ \\
\begin{CD}
A \otimes A @> \Delta \otimes \id >> A \otimes A \otimes A\\
@A{\Delta}AA @AA {\id \otimes \Delta} A\\
A  @> \Delta >> A \otimes A
\end{CD}
\end{array}
$$
A \textit{morphism} $\phi:A \lra B$ of two coalgebras,
with comultiplication  $\Delta_A$, $\Delta_B$ and counit
$\ep_A$, $\ep_B$ respectively,
is a linear map such that 
 $$(\phi \otimes \phi)  {\circ} \Delta_A
=\Delta_B {\circ} \phi,\qquad \ep_B  {\circ} \phi=\ep_A\,.$$
\item The multiplication $\mu$ and the unit $i$ are
coalgebra morphisms.
\item  The comultiplication $\Delta$ and the counit $\epsilon$
are algebra morphisms.
\item $A$ is equipped with a bijective linear map $S:A \lra A$
called the \textit{antipode} such that the following diagrams commute:
\index{antipode}
$$
\begin{CD}
A \otimes A @> {S \otimes \id} >> A \otimes A \\
@A \Delta AA   @VV \mu V \\
A @> {i {\circ} \ep} >>  A
\end{CD}
\qquad
\begin{CD}
A \otimes A @> {\id \otimes S} >> A \otimes A \\
@A \Delta AA   @VV \mu V \\
A @> {i {\circ} \ep} >>  A
\end{CD}
$$
\end{enumerate}

\smallskip

(Conditions 3 and 4 are equivalent).

\smallskip

A \textit{Hopf algebra morphism} is a linear map $\phi:A \rightarrow B$
which is a morphism of both the algebra and
 coalgebra structures of $A$
and $B$ and in addition it commutes with the antipodes
$$
S_B \circ \phi=\phi \circ S_A\,,
$$
where $S_A$ and $S_B$ denote, respectively, the antipodes in $A$ and $B$.

\smallskip

If $A$ satisfies only the first four properties is
called a \textit{bialgebra}.
\end{definition}

Bialgebras do not {necessarily} have an antipode, but if
an antipode exists it is unique.

\medskip

Let $A$ be a coalgebra with comultiplication
$\Delta:A\rightarrow A\otimes A$ and  consider the linear map
\begin{equation}
\label{eq:flip}
\begin{CD}A\otimes A@>\sigma>>A\otimes A\\a\otimes b@>>>b\otimes a\,.\end{CD}
\end{equation}
We say that $A$ is {\it cocommutative} if
$$\Delta=\sigma\circ\Delta\,.$$ In the same way, if $A$ is an algebra with multiplication $\mu:A\times A\rightarrow A$, we can express the commutativity condition as
$$\mu=\mu\circ \sigma\,.$$
We will say that a Hopf algebra is commutative or cocommutative if the underlying algebra and coalgebra structures are so.

\medskip

Let $A$ be a coalgebra with comultiplication $\Delta$ and counit $\epsilon$. We say that a subspace $I\subset A$ is a {\it coideal} \index{coideal}  if
$$\Delta(I) \subset I \otimes A+A \otimes I, \quad
\ep(I)=0\,.$$
If $A$ is a Hopf algebra, we say that $I\subset A$ is a {\it Hopf ideal} \index{Hopf ideal} if $I$ is an ideal of the algebra structure, a coideal of the coalgebra structure and
$$S(I) \subset I\,.$$
One can check immediately that in that case the algebra $A/I$ inherits naturally
a Hopf algebra structure from $A$.

\medskip

The square of the antipode of a Hopf algebra, $S^2$, is an isomorphism of Hopf algebras. However, it is not true in general that $S^2=\id$. For the cases in which the Hopf algebra is commutative or cocommutative, then we have that $S^2=\id$.

\medskip

Let $A$ be an algebra with multiplication $\mu:A\otimes A\rightarrow A$ and unit $i:k\rightarrow A$. Let us consider its dual space $A^*$. Then, the dual maps
$$\begin{CD}A^*@>\mu^*>>(A\otimes A)^*\\
a@>>>\mu^*(\alpha)=\alpha\circ\mu\,,\end{CD}\qquad \begin{CD}A^*@>i^*>> k\\
\alpha @>>>i^*(\alpha)=\alpha\circ i\,,\end{CD}$$
define on $A^*$ a coalgebra structure, provided we can identify $(A\otimes A)^*\cong A^*\otimes A^*$. If $A$ is finite dimensional, this is always the case.
The concepts of algebra and coalgebra are then seen to be dual concepts. This prompts the following definition:

\medskip

\begin{definition}\label{dualhopf} We say that the two Hopf algebras $H$ and $H'$ are
in \textit{duality} with each other if we have a non degenerate
pairing $\langle \, , \, \rangle: H \times H' \lra k$ satisfying the
properties:
\begin{eqnarray*}
    &\langle uv,x\rangle=\langle u\otimes v,\Delta'(x)\rangle,\qquad
   & \langle u,xy\rangle=\langle \Delta(u),x\otimes y\rangle,\\ \\
    &\langle 1,x\rangle=\ep'(x),\qquad
    &\langle u,1\rangle=\ep(u)\,,\\ \\
    &\langle S(u),x\rangle=\langle u,S'(x)\rangle\qquad&\text{for}\quad u,v\in H,\quad x, y\in H'\,.
    \end{eqnarray*}
$\Delta,\epsilon, S$ and $\Delta',\epsilon',S'$ denote
the comultiplication, counit and antipode in $H$ and $H'$ respectively.\hfill$\blacksquare$
\end{definition}

{
In the algebra community, it is typical to frame Pontryagin duality in terms of dually paired Hopf algebras.  In the $C^*$-algebra community, it is more common to describe it in terms of skew-paired Hopf algebras, as we have presented in Observation \ref{ex:pairing}.  The difference is that with a skew-pairing, the coalgebra structure on $H$ is replaced by its opposite via the flip map \eqref{eq:flip}, which leads us also to use the inverse of the antipode $S^{-1}$ on $H'$.  For details, see \cite[\S1.2.4 \& \S8.2.1]{KliSch}.
}

\bibliographystyle{alpha}
\bibliography{refs}

\end{document}